\documentclass[12pt]{amsart}
\usepackage{graphicx,epsfig}
\usepackage{amssymb,amsfonts,amsmath,amsthm,dsfont,wasysym,pifont,stmaryrd}
\usepackage{epstopdf,yfonts}
\usepackage{hyperref}
\newtheorem{theorem}{Theorem}[section]
\newtheorem*{theorem1'}{Theorem 1'}

\newtheorem{questions}{Question}
\newtheorem{prop}{Proposition}[section]

\newtheorem{claim}{Claim}[section]
\newtheorem*{clam}{Claim}
\newtheorem{cor}{Corollary}[section]
\newtheorem{fact}{Fact}[section]
\newtheorem{lemma}{Lemma}[section]
\newtheorem{definition}{Definition}[section]

\newcommand{\1}{\mathds{1}}
\renewcommand{\(}{\left(}
\renewcommand{\)}{\right)}

\renewcommand{\^}[1]{\widehat{#1}}
\renewcommand{\~}[1]{\overline{#1}}
\renewcommand{\geq}{\geqslant}
\renewcommand{\leq}{\leqslant}

\newcommand{\<}{\left\langle}
\renewcommand{\>}{\right\rangle}
\newcommand{\8}{\infty}

\renewcommand{\:}{\colon}
\renewcommand{\a}{\alpha}

\newcommand{\Aut}{\text{Aut}}
\renewcommand{\b}{\beta}
\newcommand{\B}{\mathcal{B}}

\newcommand{\bs}{\backslash}
\newcommand{\C}{\mathbb{C}}
\renewcommand{\Cap}[2]{\underset{#1}{\overset{#2}{\cap} }}

\renewcommand{\d}{\delta}
\newcommand{\D}{\Delta}
\newcommand{\e}{\epsilon}
\newcommand{\f}{\varphi}

\newcommand{\g}{\gamma}
\newcommand{\G}{\Gamma}

\newcommand{\GL}[2]{\textnormal{GL}_{#1}(#2)}

\renewcommand{\H}{\mathcal{H}}
\newcommand{\Hom}{\textnormal{Hom}}

\renewcommand{\int}{\varint}
\newcommand{\K}{\mathbb{K}}

\renewcommand{\l}{\left}

\newcommand{\Lim}[1]{\underset{#1}{\lim}}

\newcommand{\N}{\mathbb{N}}
\newcommand{\norm}{\trianglelefteqslant}

\newcommand{\Oplus}[2]{\underset{#1}{\overset{#2}{\oplus} }}
\newcommand{\Otimes}[1]{\underset{#1}{\otimes}}

\renewcommand{\P}{\textrm{P}}

\newcommand{\Prod}[1]{\underset{#1}{\prod}}
\newcommand{\Proj}{\textrm{Proj}}
\newcommand{\Ps}[1]{{\mathbb{P}(#1)}}

\newcommand{\PGL}[2]{\textnormal{PGL}_{#1}(#2)}

\newcommand{\Q}{\mathbb{Q}}
\renewcommand{\qedsymbol}{\fin}
\newcommand{\Quote}[1]{\textquotedblleft#1\textquotedblright}
\renewcommand{\r}{\right}
\newcommand{\R}{\mathcal{R}}
\renewcommand{\Re}{\mathbb{R}}

\newcommand{\s}{\sigma}
\renewcommand{\S}{\Sigma}
\newcommand{\set}[2]{\left \{ #1 \colon #2 \right \}}
\newcommand{\sdp}[2]{#1 \ltimes #2}
\newcommand{\SL}[2]{\textnormal{SL}_{#1}(#2)}

\newcommand{\stab}{\mathrm{stab}}
\newcommand{\Sum}[2]{\underset{#1}{\overset{#2}{\sum} }}
\newcommand{\Sup}[1]{\underset{#1}{\sup}}
\newcommand{\T}{\mbox{$\mathbb{T}$}}
\renewcommand{\t}{\tau}
\newcommand{\tr}{\textsl{tr}}
\newcommand{\Un}{\mbox{$\mathcal{U}$}}

\newcommand{\Z}{\mathbb{Z}}

\title{Relative Property (T) and Linear Groups}
\author{Talia Fern\'os}
\address{University of Illinois at Chicago \\
Dept. of MSCS (m/c 249) \\
 851 South Morgan Street \\ 
 Chicago, IL 60607-7045}
\email{talia@math.uic.edu}
\begin{document}
\maketitle

\begin{abstract}

Relative property (T) has recently been used to construct a variety of new rigidity phenomena, for example in von Neumann algebras and the study of orbit-equivalence relations. However, until recently there were few examples of group pairs with relative property (T) available through the literature. This motivated the following result: A finitely generated group $\G$ admits a $\Re$-special linear representation with non-amenable $\Re$-Zariski closure if and only if it acts on an Abelian group $A$ (of finite nonzero $\Q$-rank) so that the corresponding group pair $(\Gamma \ltimes A, A)$ has relative property (T).

The proof is constructive. The main ingredients are Furstenberg's celebrated lemma about invariant measures on projective spaces and the spectral theorem for the decomposition of unitary representations of Abelian groups. Methods from algebraic group theory, such as the restriction of scalars functor, are also employed. 

\end{abstract}


\section{Introduction}

Recall that if $\G$ is a group and $A \leq \G$ is a closed subgroup then the group pair $(\G, A)$ is said to have relative property (T) if every unitary representation of $\G$ with almost invariant vectors has $A$-invariant vectors. And $\G$ is said to have property (T) if $(\G, \G)$ has relative property (T)%
\footnote{We will assume throughout this paper that groups are locally compact and second countable, Hilbert spaces are separable, unitary representations are strongly continuous (in the usual sense), fields are of characteristic 0, and local fields are not discrete. Furthermore, all countable groups will be given the discrete topology, unless otherwise specified. }.

In 1967 D. Kazhdan used the relative property (T) of the group pair $(\sdp{\SL{2}{\K}}{\K^2}, \K^2)$ to show that $\SL{3}{\K}$ has property (T), for any local field $\K$ \cite[Lemmas 2 \& 3]{Kazhdan}. Later in 1973 G. A. Margulis used the relative property (T) of $(\sdp{\SL{2}{\Z}}{\Z^2}, \Z^2)$ \cite[Lemma 3.18]{Margulis} in order to construct the first explicit examples of families of expander graphs. It was he who later coined the term.

Recently relative property (T) has been used to construct a variety of new phenomena. Most notable is the recent work of S.Popa. He has shown that every countable subgroup of $\Re_+^*$ is the fundamental group of some II$_1$-factor \cite{Popa3}, and constructed examples of II$_1$ factors with rigid Cartan subalgebra inclusion \cite{Popa1}. Also D. Gaboriau with S. Popa constructed uncountably many non-orbit equivalent (free and ergodic measure-preserving) actions of the free group $F_n$ (for $n \geq 2$) on the standard probability space. See \cite{PopaGab} and \cite{Popa2} and the references contained therein. 

In a completely different direction, A. Navas, extending his previous work with property (T) groups, showed that relative property (T) group pairs acting on the circle by $C^2$ diffeomorphisms are trivial, in a suitable sense \cite{Navas}. Also, M. Kassabov and N. Nikolov \cite[Theorem 3]{KasNik} used relative property (T) to show that $\SL{n}{\Z[x_1, \dots, x_k]}$ has property ($\tau$) for $n \geq 3$.

We also refer to A. Valette's paper \cite{Valette} for more applications concerning, for example, the Baum-Connes conjecture. 

Unfortunately, until recently the examples of group pairs with relative property (T) available in the literature have been scarce:

\begin{itemize}
  \item If $n \geq 2$ then $(\sdp{\SL{n}{\Re}}{\Re^n}, \Re^n)$ and $(\sdp{\SL{n}{\Z}}{\Z^n}, \Z^n)$ have relative property (T). \cite[10-Proposition]{HarpeValette}
  \item If $\G \leq \SL{2}{\Z}$ is not virtually cyclic then $(\sdp{\G}{\Z^2}, \Z^2)$ has relative property (T). \cite[Example 2 Section 5]{Burger}
  \item And, now, in a recent preprint of A. Valette \cite{Valette}: If $\G$ is an arithmetic lattice in an absolutely simple Lie group then there exists a homomorphism $\G \to \SL{N}{\Z}$ such that the corresponding pair $(\sdp{\G}{\Z^N}, \Z^N)$ has relative property (T). 
 \end{itemize}

We remark that $\sdp{\SL{n}{\Re}}{\Re^n}$ actually has property (T) for $n \geq 3$ \cite{Wang} and so $(\sdp{\SL{n}{\Re}}{\Re^n}, A)$ has relative property (T) for any closed $A \leq \sdp{\SL{n}{\Re}}{\Re^n}$. Indeed, if $A \leq G \leq H$ are groups, and $G$ has property (T) then $(H,A)$ has relative property (T). 

On the other extreme, if $S$ is an amenable group then $(S,A)$ has relative property (T) if and only if $A$ is compact. (See Lemma 8.3 in Section 8.) So, if one wants to find new examples of group pairs with relative property (T), they should not rely on the property (T) on one of the groups in question and they should be of the form $(\G, A)$ where $\G$ is non-amenable and $A$ is amenable but not compact. 

Using these examples as a guide, one may ask to what extent can group pairs with relative property (T) be constructed? We offer the following as an answer to this question:

\begin{theorem}
Let $\G$ be a finitely generated group. The following are equivalent:
\begin{enumerate}
  \item There exists a homomorphism $\f : \G \to \SL{n}{\Re}$ such that the $\Re$-Zariski-closure $\~{\f(\G)}^Z(\Re)$ is non-amenable.
  \item There exists an Abelian group $A$ of nonzero finite $\Q$-rank and a homomorphism $\f' \: \G \to \Aut(A)$ such that the corresponding group pair $(\sdp{\G}{_{\f'}A}, A)$ has relative property (T).
\end{enumerate}
\end{theorem}

\textbf{Remark:} In the direction of $(1) \implies (2)$, more information can be given. Namely, we will specifically find that $A = \Z[S^{-1}]^N$ where $S$ is some finite set of rational primes, as is pointed out below. Also in the direction of $(2) \implies (1)$ we will find that $A$ can be taken to be of the form $\Z[S^{-1}]^N$.  

\break
\subsection{Outline of the proof of Theorem 1 in the direction $(1) \implies (2)$}

\subsubsection{From Transcendental to Arithmetic}

This step is a matter of showing that from an arbitrary representation $\f \: \G\to \SL{n}{\Re}$, such that the $\Re$-Zariski closure $\~{\f(\G)}^Z(\Re)$ is non-amenable, we may find an arithmetic representation $\psi \: \G \to \SL{m}{\Q}$ such that the $\Re$-Zariski closure $\~{\psi(\G)}^Z(\Re)$ is non-amenable.
  
\subsubsection{Relative Property (T) for $\Re^N$.}

We establish the existence of a subgroup $\G_0 \norm \G$ of finite index and a ``nice'' representation $\a \: \G_0 \to \SL{N}{\Q}$ such that $(\sdp{\G_0}{_{\a} \Re^N}, \Re^N)$ has relative property (T). The representation $\a$ is a factor of $\psi|_{\G_0}$. 

\subsubsection{Fixing the Primes}

  We show that, after conjugating the representation $\a$ by an element in $\GL{N}{\Q}$ if necessary, we may assume that $\a \: \G_0 \to \SL{N}{\Z[S^{-1}]}$ and that $\a(\G_0)$ is not $\Q_p$-precompact for each $p \in S$. The representation $\a$ is so nice that this allows us to conclude that $(\sdp{\G_0}{_{\a} \Q_p^N}, \Q_p^N)$ has relative property (T) for each $p \in S$.

\subsubsection{Products and Induction} 
  
  The set $S$ of primes in Step 3 is finite, and we show that the relative property (T) passes to finite products. Namely, if $(\sdp{\G_0}{_{\a} \Q_p^N}, \Q_p^N)$ has relative property (T) for each $p \in S\cup \{\8\}$ then setting $V = \Prod{ p \in S\cup \{\8\} }\Q_p^N$ we have that $\( \sdp{\G_0}{V }, V \)$ has relative property (T). 

Let $A = \Z[S^{-1}]^N$ and recall that the diagonal embedding $A \subset V$ is a lattice embedding. Since $\a(\G_0) \leq \SL{N}{\Z[S^{-1}]}$ we have that $\G_0$ acts on $A$ by automorphisms. Since $\sdp{\G_0}{A}$ is a lattice in $\sdp{\G_0}{V}$ we have that $(\sdp{\G_0}{A}, A)$ has relative property (T).

\subsubsection{Extending up from a finite index subgroup}
  
  We show that if $k = [\G : \G_0]$ then there is a homomorphism  $\a' \: \G \to \SL{kN}{\Z[S^{-1}]}$ such that $( \sdp{\G}{A^k}, A^k)$ has relative property (T). 

\subsection{Outline of the proof of Theorem 1 in the direction $(2) \implies (1)$}

\subsubsection{Managing $A$}
We choose $A$ to be of minimal (non-zero) $\Q$-rank among all Abelian groups satisfying condition (2). Under the hypothesis, we show that we may assume that $A$ is torsion free and hence a subgroup of $\Q^n$ where $n$ is the $\Q$-rank of $A$. This yields that there are finite sets of primes $S_i$ such that, up to isomorphism, $A = \Oplus{i = 1}{n} \Z[S_i^{-1}]$.

\subsubsection{An Invariant subgroup of $A$}
We choose $m \in \{1, \dots, n\}$ such that $|S_m| \geq |S_i|$ for each $i \in \{1, \dots, n\}$. Letting $I_m = \set{i}{S_i = S_m}$ we get that $A_m = \Oplus{i \in I_m}{} \Z[S_m^{-1}]$ is $\G$-invariant. By minimality of $A$ it follows that $A = A_m \cong \Z[S_m^{-1}]^n$. Set $S = S_m$.

\subsubsection{A is a lattice}
Let $V =  \Re^n \times \Prod{p \in S}\Q_p^n$. Since $A \subset V$ is a co-compact lattice it follows that $(\sdp{\G}{V}, V)$ has relative property (T). 

\subsubsection{The $\Re$-component}
Since $\Prod{p \in S}\Q_p^n \subset V$ is $\G$-invariant we have that $(\sdp{\G}{\Re^n},\Re^n)$ has relative property (T). 

\subsubsection{The Image of $\G$}
If $\f \: \G \to \GL{n}{\Q}$ is the corresponding homomorphism, then $\ker(\f) \norm \sdp{\G}{\Re^n}$ so that $(\sdp{\f(\G)}{\Re^n},\Re^n)$ has relative property (T). 

\subsubsection{The Zariski Closure}
If $(\sdp{\f(\G)}{\Re^n},\Re^n)$ has relative property (T) then $(\sdp{\~{\f(\G)}^Z(\Re)}{\Re^n},\Re^n)$ has relative property (T). It is shown that this implies that $\~{\f(\G)}^Z(\Re)$ is not amenable. 

\subsection{Organization of the Paper}

We present the paper in the following order:

\subsubsection{Section 2}
In Section 2 we discuss some algebraic preliminaries in order to make the rest of the exposition consistent and coherent. 

\subsubsection{Section 3}
In Section 3 we state and discuss the main theorems (Theorem 2 and Theorem 3) that will be used in the proof of Theorem 1 in the direction of $(1) \implies (2)$. Their roles are:
\begin{itemize}
  \item []  \begin{itemize}
  \item [\textit{ Thm. 2}] To give a criterion on a group $\G$ (we will call it Property (F$_p$)) for which we may construct group pairs $(\sdp{\G}{\Q_p^n}, \Q_p^n)$ having relative property (T).
  \item [\textit{ Thm. 3}] To give a criterion on a group $\G$ for which there is a finite set of primes $S$ such that we may construct group pairs $(\sdp{\G}{\Z[S^{-1}]^n}, \Z[S^{-1}]^n)$ having relative property (T).
\end{itemize}
\end{itemize}

\subsubsection{Section 4}
In Section 4 we prove Theorem 2.

\subsubsection{Section 5}
In Section 5 we prove Theorem 3 using Theorem 2. 

\subsubsection{Section 6}
In Section 6 we prove an algebro-geometric specialization proposition (Proposition 4). It exactly yields step 1 in the proof of Theorem 1 for the direction $(1) \implies (2)$.

\subsubsection{Section 7}
 In Section 7 we prove Theorem 1 in the direction of $(1) \implies (2)$ essentially as a consequence of Proposition 6.1 and Theorem 3. 

\subsubsection{Section 8}
In Section 8, we prove Theorem 1 in the direction of $(2) \implies (1)$. The proof is simple, and is pretty much self contained. 

\subsection{Acknowledgments:} I'd like to thank Alex Furman for being a truly excellent advisor. In particular he deserves a great deal of thanks for his many detailed readings of this paper and his instructive comments and suggestions. He also proposed the original idea behind this work. I'd also like to thank Alain Valette for sending me a preprint of his paper \cite{Valette}. It came at an opportune time as it allowed for the generalization of the work I had in progress. I'd also like to thank him for his comments on this work.

This work is a part of my doctoral thesis.

\section{Algebraic Preliminaries}

\subsection{A word about Zariski Closures}\cite[Section AG.13]{Borel}, \cite[Section3.1]{Zimmer}

Let $k$ be a field and $\~K$ an algebraically closed field containing $k$. Recall that to every subset $V \subset \~K^n$ there corresponds an ideal $I_{\~K}(V) \subset \~K[x_1, \dots, x_n]$ such that $p \in I_{\~K}(V)$ if and only if $p|_V \equiv 0$. The set $V$ is said to be Zariski closed if $V = \set{a \in \~K^n}{p(a) = 0 \text{ for every } p \in I_{\~K}(V)}$, that is, if it is exactly the zero-set of its ideal. 

Furthermore, $V$ is said to be defined over $k$ if there exists an ideal $I_k(V) \subset k[x_1, \dots, x_n]$ such that $I_k(V) \cdot \~K[x_1, \dots, x_n] = I_{\~K}(V)$. In such a case we write  \begin{equation}
V(k) := \set{a \in k^n}{p(a) = 0 \text{ for every }p \in I_k(V)} \nonumber
\end{equation}
to denote the $k$-points of $V$. Observe that it could happen that $V(k) = \phi$ despite the fact that $I_k(V) \neq k[x_1, \dots, x_n]$. (Take for example $\~K= \C$ and $k = \Re$ and $V = \{ i, -i\} \subset \~K^1$. Then $I_\Re(V) = (x^2 +1)$ is defined over $\Q$ and $V(\Re) = \phi$. This is why we need to work with algebraically closed fields to begin with!) Fortunately, the situation for groups is significantly better. 

Recall that $\GL{n}{\~K}$ is an algebraic (i.e. Zariski closed) group defined over $\Q$. 

\begin{prop}\cite[Proposition 3.1.8]{Zimmer} 
Suppose that $G(\~K) \leq \GL{n}{\~K}$ is an algebraic group such that $G(k) := \GL{n}{k}\cap G(\~K)$ is Zariski dense in $G(\~K)$. Then $G(\~K)$ is defined over $k$.  
\end{prop}

\begin{prop}[Chevalley]\cite[Theorem 3.1.9]{Zimmer} If $G(\~K)$ is an algebraic group defined over $k$ then $G(k)$ is Zariski dense in $G(\~K)$.
\end{prop}

Note that this means in particular, that if $G(\~K) \leq \GL{n}{\~K}$ is Zariski closed, nontrivial, and defined over $k$ then $G(k)$ is nontrivial as well!

Now, if $\G \leq \GL{n}{k}$ is any subgroup, then the $\~K$-Zariski closure is denoted by $\~{\G}^Z(\~K)$. We say $\~K$-Zariski closure since this depends on the algebraically closed field $\~K$. Indeed, if $\~K'$ is another algebraically closed field containing $k$, then by the above propositions, $\G$ is also Zariski dense in $\~{\G}^Z(\~K')$. 

Observe that this notion is well defined even if the field is not algebraically closed. Namely, let $F$ be a field containing $k$ and let $\~F$ be its algebraic closure. We define the $F$-Zariski closure of $\G$ to be $\~{\G}^Z(F) := \~{\G}^Z(\~F) \cap \GL{n}{F}$. In general we make use of this when it has additional topological content. For example if $k= \Q$ and $F = \Q_p$ for some prime $p$. Then the group $\~{\G}^Z(F)$ is a $p$-adic group and has a lot of nice additional structure.

 
 
\subsection{Restriction of Scalars:}

Let $K$ be a finite separable extension of a field $k$ (of any characteristic) and $\S := \{ \s \: K \to \~k\}$ be the set of $k$-linear embeddings of $K$ into $\~k$ a fixed separable closure of $k$. There is a functor called the \emph{restriction of scalars} functor which maps the category of linear algebraic $K$-groups and $K$-morphisms into the category of linear algebraic $k$-groups and $k$-morphisms. Namely, let $H$ be an algebraic $K$-group defined by the ideal $I \subset K[X]$. Then, for each $\s \in \S$ the algebraic group $^\s H$ is defined by $\s(I) \subset \s(K)[X]$, the ideal obtained by applying $\s$ to the coefficients of the polynomials in $I$. The restriction of scalars of $H$ is $\R_{K/k}H \cong \Prod{\s \in G}^\s H$. It has the following properties \cite[Section 6.17]{BorelTits} \cite[ Proposition 6.1.3] {Zimmer} \cite[Section 12.4 ]{Springer}:
\begin{enumerate}
  \item There is a $K$-morphism $\a \: \R_{K/k}H \to H$ such that the pair $(\R_{K/k}H , \a)$ is unique up to $k$-isomorphism.
  \item If $H'$ is a $k$-group and $\b \: H' \to H$ is a $K$-morphism then there exists a unique $k$-morphism $\b' \: H' \to \R_{K/k}H $ such that $\b = \a\circ \b'$.
\item  If $K'$ is any field containing $K$ then $\R_{K/k}H (K') \cong \Prod{\s \in G}^\s H(K')$.
 \item The algebraic type of the group is respected. Namely, if $H$ has the property of being reductive (respectively semi-simple, parabolic, or Cartan) then $\R_{K/k}H$ is reductive (respectively semi-simple, parabolic, or Cartan).
 \item The algebraic type of subgroups is respected. Namely, if $P  \leq H$ is a $K$-Cartan subgroup (respectively $K$-maximal torus, $K$-parabolic subgroup) then $\R_{K/k}P \leq \R_{K/k}H$ is a $k$-Cartan subgroup (respectively $k$-maximal torus, $k$-parabolic subgroup).
  \item There is a correspondence of rational points: Consider the diagonal embedding $\D \:  H(K) \to \Prod{\s \in \S} ^\s H(K)$ defined pointwise by $h \mapsto \Prod{\s \in \S} \s(h)$. Then we have the correspondence  $\R_{K/k}H (k) \cong \D(H(K))$. 

\end{enumerate}

\noindent
\textbf{ Disclaimer:} In the sequel we consider the isomorphism $\R_{K/k}H \cong \Prod{\s \in \S}^\s H$ as equality.


\section{The Main Theorems 2 and 3}


Note that if $\G$ is a finitely generated group and $\f \: \G \to \SL{n}{\~\Q}$ is an algebraic representation, then there is a field $K_\f$ which is a normal finite extension of $\Q$ such that $\f(\G) \leq \SL{n}{K_\f}$. (Take for example, the normal field generated by the entries of some finite generating set for $\f(\G)$.)

With this notation in place, we give the following definition, which will be used to find group pairs with relative property (T). 

\begin{definition}
Let $\G$ be a finitely generated non-amenable group and $p \in \{2, 3, 5, \dots, \8\}$ a rational prime. Then $\G$ is said to satisfy property (F$_p$) (after Furstenberg) if there exists an algebraic homomorphism $\f \: \G \to \SL{n}{\~\Q}$ satisfying the following conditions: %
\begin{enumerate}
  \item The $\~\Q$-Zariski closure $H = \~{\f(\G)}^Z(\~\Q)$ is $\~\Q$-simple.
  \item There are no $\f(\G)$-fixed vectors.
  \item The natural diagonal embedding $\D: \f(\G) \to \R_{K_\f/\Q}H(\Q)$ is not pre-compact in the $p$-adic topology.
\end{enumerate} 

In such a case, we say that the representation $\f$ realizes property (F$_p$) for $\G$

\end{definition}

Recall that the archimedean valuation on $\Q$ is called the prime at infinity. So, according to convenience, we use both notations $\Re$ and $\Q_\8$ to denote the completion of $\Q$ with respect to the archimedean valuation. 
  
 
 \begin{theorem}
Let $\G$ be a group satisfying property (F$_p$). Then, there exists a rational representation $\f' \: \G \to \SL{N}{\Q}$ such that $(\sdp{\G}{_{\f'} \Q_p^N}, \Q_p^N)$ has relative property (T).
\end{theorem}


 \begin{theorem}
Suppose that $\G$ is a group with property (F$_\8$). Then there exists a finite set of primes $S \subset \Z$ and a representation $\rho \: \G \to \SL{N}{\Z[S^{-1}]}$ such that, if $A = \Z[S^{-1}]^N$ then $(\sdp{\G}{_{\rho}A}, A)$ has relative property (T). 
\end{theorem}

 
 \textbf{Remark:} Conditions (1) and (2) of property (F$_p$) can be seen as an irreducibility requirement. With this in mind, we see that Theorems 2 and 3 say that irreducibility and unboundedness are sufficient ingredients to cook up a relative property (T) group pair.


\section{Theorem 2}


\subsection{How to Find Relative Property (T)}

 Our first task is to establish a sufficient condition for the presence of relative property (T); one that lends itself to the present context. The following is due to M. Burger [Propositions 2 and 7]\cite{Burger}. In what follows $\K$ is a local field and $\^\K \cong \Hom(\K, S^1)$ is the unitary dual. Recall that $\^\K$ is topologically isomorphic to $\K$ \cite[Theorem 7-1-10 ]{Goldstein}. As such we will often not distinguish between $\GL{n}{\K}$ and $\GL{n}{\^\K}$.



 \begin{prop}[Burger's Criterion for Relative Property (T)]
Suppose that $\f \: \G \to \GL{N}{\K}$ is such that there is no $\G$-invariant probability measures on $\Ps{\^\K^N}$. Then, $(\sdp{\G}{_\f \K^N}, \K^N)$ has relative property (T).
\end{prop}

\begin{proof}

Let $\rho \: \sdp{\G}{\K^N} \to \Un(\H)$ be a unitary representation with $\G$-almost invariant vectors and $\P \: \B(\^\K^N) \to \Proj(\H)$ the projection valued measure associated to $\rho|_{\K^N}$, where $\B(\^\K^N)$ denotes the Borel $\s$-algebra of $\^\K^N$. Recall that $\P$ has the following properties:

\begin{enumerate}
\item $\P(\^\K^N) = \mathrm{Id}$
\item For every $v \in \H$ the measure $B \mapsto \<\P(B)v, v\>$ is a positive Borel measure on $\^\K^N$ with total mass $\|v\|^2$. 
  \item For every $\g \in \G$ we have that \begin{center}
$\rho(\g^{-1}) \P(B) \rho(\g) = \P(\g^*B)$.
\end{center}
  \item The projection onto the subspace of $\K^N$-invariant vectors is $\P(\{0\})$. 
\end{enumerate}

Let $v_n \in \H$ be a sequence of $(\e_n, F_\G)$-almost invariant unit vectors where $\e_n \to 0$ and $F_\G$ is a finite generating set for $\G$. Define the probability measures $\mu_n(B):= \<\P(B)v_n, v_n\>$.

\begin{clam}
The sequence of measures $\{\mu_n\}$ is almost $\G$-invariant. Namely $\| \g_*\mu_n - \mu_n\| \leq 2\e_n$ for each $\g \in F_\G$.
\end{clam}

\begin{proof}

Let $B \subseteq \^\K^N$ be a Borel set and $\g \in F_\G$. Then 
\begin{equation}
\begin{split}
|\mu_n(\g^* B) - \mu_n(B)| & =  |\< \pi(\g^{-1})\P(B)\pi(\g)v_n, v_n \> - \<\P(B)v_n, v_n \>|  \\ 
& \leq  | \< \pi(\g^{-1}) \P(B) \pi(\g) v_n, v_n \> - \< \pi(\g^{-1}) \P(B) v_n,  v_n \> |  \\
& \qquad \qquad \qquad \qquad \qquad \qquad \quad + | \< \pi(\g^{-1}) \P(B) v_n, v_n \> - \<  \P(B) v_n, v_n \>|  \\  
&= |\<\pi(\g^{-1}) \P(B) (\pi(\g)v_n - v_n), v_n\>| + |\<\P(B)v_n, (\pi(\g)v_n - v_n)\>|  \\ 
&\leq  \| \pi(\g^{-1})\P(B)\| \cdot \| \pi(\g)v_n - v_n \| + \|\P(B)\| \cdot \| \pi(\g)v_n - v_n \| \leq 2\e_n. 
\end{split} \nonumber  
\end{equation} 

Thus the sequence of probability measures $\{\mu_n\}$ is almost $\G$-invariant. 
\end{proof}

Suppose by contradiction that the group pair $(\sdp{\G}{\K^N},\K^N)$ fails to have relative property (T). Then for each $n$, $\mu_n(\{0\}) = 0$. This allows us to pass to the associated projective space. 

Namely let $p \: \^\K^N \bs \{0\} \to \Ps{\^\K^N}$ be the natural projection. Define the probability measures $\nu_n := p_*\mu_n$.
It is clear that they also satisfy the following inequality for any $\g \in F_\G$:

\begin{equation}
\| \g_*\nu_n - \nu_n\| \leq 2\e_n
\nonumber
\end{equation}

Exploiting the compactness of $\Ps{\^\K^N}$, we get that a weak-$*$ limit point of $\{\nu_n\}$ will necessarily be $\G$-invariant, a contradiction of the hypothesis that there are no $\G$-invariant probability measures on $\Ps{\^\K^N}$. 
\end{proof}

This is a powerful criterion when taken together with the following:

\begin{lemma}[Furstenberg's Lemma] \cite[Lemma 3.2.1, Corollary 3.2.2]{Zimmer}
Let $\mu$ be a Borel probability measure on $\Ps{\K^N}$. Suppose that $\G \leq \PGL{N}{\K}$ leaves $\mu$ invariant. If $\G$ is not precompact then there exists a nonzero subspace $V \subsetneq \K^N$ which is invariant under a finite index subgroup of $\G$ and such that $\mu[V] >0$. 
\end{lemma}

These two statements will be used to show the presence of relative property (T) once we have a nice representation to work with. The representation will be provided by the following considerations.



\subsection{The Tensor Representation:}

Let $K$ be a finite normal extension of $\Q$ with Galois group $G$. Consider the vector space $W(K) =\underset{\s \in G}{\otimes} K^n$ and the representation of $\R_{K/\Q}\SL{n}{K} \cong \Prod{ \s \in G} ^\s \SL{n}{K}$ on $W(K)$, defined by $\t : \Prod{\s \in G} g_\s \mapsto \Otimes{\s \in G} g_\s$. This induces a representation $\D_\t \: \SL{n}{K} \to \SL{}{W(K)}$ defined by $\D_\t = \t \circ \D$.

There are two reasons which make this an excellent representation to work with. The first is due to Y. Benoist and is taken from \cite[Lemma 1]{Valette}. 

\begin{lemma}
The faithful representation $\D_\t \: \SL{n}{K} \to \SL{}{W(K)}$ is defined over $\Q$ and there is a $\Q$-subspace $W(\Q)$ of $W(K)$ such that the map $K \otimes W(\Q) \to W(K)$ is an $\SL{n}{K}$-equivariant isomorphism. 
\end{lemma}

The second reason is observed in \cite[Item 1, page 9]{Valette}:

\begin{lemma}
If $H(K) \leq \SL{n}{K}$ is a group without fixed vectors in $K^n$ then for each $\s_0 \in G$ the restricted representation $\tau_0 = \tau \l |_ {{^{\s_0}} {H(K)}} \r.  \: ^{\s_0} H(K) \to \SL{}{W(K)}$ also has no invariant vectors.
\end{lemma}

\begin{proof}

Although we are thinking of $^{\s_0}H(K)$ as being a subgroup of $\SL{n}{K}$, for the sake of clarity it is necessary to denote by $\rho_0 \: {^{\s_0}H(K)} \to \SL{n}{K}$ the identity representation, so that $\rho_0(^{\s_0}H(K)) = {^{\s_0}H(K)}$. 

With this notation, it is clear that $\tau_0 \: ^{\s_0}H(K) \to \SL{}{W(K)}$ is given by $\tau_0 = \rho_0 \Otimes{\s \neq \s_0} \1$, where $\1$ denotes the trivial representation. Namely, $^{\s_0}H(K)$ acts trivially on each tensor-factor except the one corresponding to $\s_0$, where it acts via $\rho_0$.

Also recall the fact that 
\begin{equation}
\Otimes{\s \in G} K^n \cong \(\Otimes{\s \neq \s_0} K^n\) \otimes K^n \cong \Hom \( \( \Otimes{\s \neq \s_0} K^n \)^*, K^n \). \nonumber
\end{equation} 

Under this isomorphism, a vector which is $^{\s_0}H(K)$-invariant corresponds to a $K$-linear map which intertwines $( \Otimes{\s \neq \s_0} \1^*, ( \Otimes{\s \neq \s_0} K^n )^* )$ with $(\tau_0, K^n)$. Since the dual of a trivial representation is trivial, it follows that the image of such a map consists of $\rho_0( ^{\s_0}H)$-invariant vectors. 

We then have that $(\tau_0, \Otimes{\s \in G} K^n)$ contains a \emph{non-zero} $^{\s_0}H(K)$-invariant vector if and only if $(\rho_0, K^n)$ contains the trivial representation; that is, if and only if $(\rho_0, K^n)$ contains a $^{\s_0}H(K)$-invariant vector. And, since $H(K)$ does not have invariant vectors in $K^n$ neither does $^{\s_0}H(K)$.
\end{proof}

Before the proof of Theorem 2, we establish a little more notation: Let $F$ be a field containing $\Q$. Then we write $W(F) = W(\Q) \otimes F$. If $F$ contains $K$ then naturally $W(F) \cong \underset{\s \in G}{\otimes} F^n$. 


\subsection{The Proof of Theorem 2:}

We retain the notation established above. Recall that if $\G$ is a group satisfying property (F$_p$) then there is a field $K$ which is a finite normal extension of $\Q$  and a representation $\f \: \G \to \SL{n}{K}$ such that %
\begin{enumerate}
  \item The Zariski-closure $H = \~{\f(\G)}^Z$ is $\~\Q$-simple.
  \item There are no $\f(\G)$-fixed vectors.
  \item The natural diagonal embedding $ \D \: \f(\G) \to \R_{K/\Q}H(\Q)$ is not pre-compact in the $p$-adic topology.
\end{enumerate} 

\begin{proof}
Consider the representation of $\f' \: \G \to \SL{}{W(\Q)}$ which is defined as $\f' = \tau \circ \D \circ \f$. We claim that $(\sdp{\G}{_{\f'} W(\Q_p)}, W(\Q_p))$ has relative property (T). 

If not then by Burger's Criterion (Proposition 4.1) there exists a $\G$-invariant probability measure $\mu$ on $\Ps{W(\^{\Q_p})}$. Since $\f'$ factors through the diagonal embedding in item (3) above, it follows that $\f'(\G) \leq \SL{}{W(\Q_p)}$ is not pre-compact, and hence the corresponding projective image in $\PGL{}{W(\^{\Q_p})}$ is also not pre-compact (since $\SL{}{W(\^{\Q_p})}$ has finite center). By Furstenberg's Lemma, there exists a non-trivial subspace $V \subsetneq W(\^{\Q_p})$ such that %
\begin{enumerate}
  \item There is a subgroup of finite index in $\G$ which preserves $V$.
  \item The mass $\mu[V]>0$.
  \item $V$ is of minimal dimension among all subspaces satisfying (1) and (2).
\end{enumerate}
We aim to show that this is impossible:

Observe that $V$ is actually $\R_{K/\Q}H(\Q_p)$-invariant. Indeed, since preserving a subspace is a Zariski-closed condition (consider the corresponding parabolic subgroup), if $\G$ has a finite index subgroup which preserves $V$ then so must the Zariski-closure $\R_{K/\Q}H(\Q_p)$. Since $H$ is $\~\Q$-simple, it is Zariski-connected and therefore so is $\R_{K/\Q}H(\Q_p)$. It follows that all of $\R_{K/\Q}H(\Q_p)$, and in particular $\G$, preserves $V$.

We claim that the map $\R_{K/\Q}H(\Q_p) \to \SL{}{V}$ is a faithful continuous homomorphism. Continuity is automatic because the representation is linear. (Observe that the semisimplicity of $\R_{K/\Q}H(\Q_p)$ guarantees that the image is in $\SL{}{V}$ versus $\GL{}{V}$.)

Since $\f'(\G) \leq \SL{}{W(\Q)}$ it follows that the subspace $V$ is defined over an algebraic field $F \subset \~\Q$, and we may as well assume that  $K \subset F$. Let $V(F)$ be the $F$-span of an $F$-basis of $V$. Then, we have the representation $\R_{K/\Q}H(F) \to \SL{}{V(F)}$. 

Recall that property (3) of the restriction of scalars says that $\R_{K/\Q}H(F) \cong \Prod{\s \in G} ^\s H(F)$, where $G$ is the Galois group of $K/\Q$. Now observe that since each $^\s H$ is $\~\Q$-simple, the kernel is either trivial, or contains  $^{\s_0} H(F)$ for some $\s_0 \in G$. Assume that the kernel is not trivial. This means that $^{\s_0}H(F)$ acts trivially on $V(F)$, i.e. that each vector in $V(F)$ is fixed by $^{\s_0}H(F)$. We claim that this is impossible:

Indeed, by Lemma 4.3, there are no $^{\s_0}H(K)$-invariant vectors in $W(K)$.  This means that $W(F)$ cannot have $^{\s_0}H(F)$-invariant vectors. If $v \in W(F)$ is $^{\s_0}H(F)$-invariant then it is $^{\s_0}H(K)$-invariant which means that $v \in W(K)$ (since the equations for $v$ are linear with coefficients in $K$), a contradiction.

Thus, the representation $\R_{K/\Q}H(\Q_p) \to \SL{}{V}$ is faithful and continuous. Since $\D \circ \f (\G) \leq \R_{K/\Q}H(\Q_p) $ is not precompact, it follows that the corresponding representation $\G \to \SL{}{V}$ is also not precompact. 

Now, consider the induced measure:
\begin{equation}
\mu_0(B) = \mu(B \cap [V]) / \mu[V]. \nonumber
\end{equation}

It is clearly $\G$-invariant. Furthermore, since $V$ was chosen to be of minimal dimension by Furstenberg's lemma,  it follows that the image of $\G$ in $\PGL{}{V}$ is pre-compact, which is a contradiction.

Thus, there are no $\G$-invariant probability measures on $\Ps{W(\^{\Q_p})}$ and so by Burger's Criterion, the group pair $(\sdp{\G}{W(\Q_p)}, {W(\Q_p)})$ has relative property (T). 

\end{proof}


\section{Theorem 3}

Recall that if $\G$ has property (F$_\8$) then there exists a representation $\f \: \G \to \SL{n}{K}$ (with $d = [K:\Q] <\8$) such that
\begin{enumerate}
  \item The Zariski-closure $H = \~{\f(\G)}^Z$ is $\~\Q$-simple.
  \item The representation $\f$ does not contain the trivial representation, that is, there are no $\f(\G)$-fixed vectors.
  \item[*(3)] The natural diagonal embedding $\D: \f(\G) \to \R_{K/\Q}H(\Q)$ is not pre-compact in the $\8$-adic (that is the $\Re$) topology.
\end{enumerate} 
\break
We now turn to the proof of Theorem 3:

\begin{proof}

Let $N = n^d$. We retain the notation from the proof of Theorem 2 and set $\Q^N \cong W(\Q)$. Recall that this gives rise to:
\begin{equation}
\f' \: \G \overset{\f}{\to} H(K) \overset{\D}{\hookrightarrow} \R_{K/\Q}H(\Q) \overset{\tau}{\hookrightarrow} \SL{N}{\Q}
\nonumber
\end{equation}
and $(\sdp{\G}{_{\f'}\Re^N}, \Re^N)$ has relative property (T) by Theorem 2.

Note that proof of Theorem 2 also shows that if there exists a prime $p$ such that condition *(3) holds at $p$ (that is if $\D \circ\f(\G)$ is also not precompact in the $p$-adic topology) then $(\sdp{\G}{_{\f'}\Q_p^N}, \Q_p^N)$ has relative property (T). (For the same $\f'$!) Let $S \subset \Z$ be the set of primes such that if $p \in S$ then condition *(3) holds at $p$. 

Next, let $S_0 \subset \Z$ be the set of primes such that if $p \in S_0$ then $p$ appears as a denominator in some entry of $\f'(\G)$. Since $\G$ is finitely generated, $S_0$ is finite and by definition $\f'(\G) \leq \SL{N}{\Z[S_0^{-1}]}$. 

Recall that going to infinity in the $p$-adic topology amounts to being \Quote{increasingly divided by $p$}. By observing that $\tau$ is faithful, we see that $S \subset S_0$ and so $S$ is also finite. Consider the following:


 \begin{lemma}
Let $S$ and $S_0 = S \cup \{p\}$ be two distinct sets of primes. If $\G \leq \SL{N}{\Z[S_0^{-1}] }$ is such that the natural embedding $\G \leq \SL{n}{\Q_p}$ is precompact, then there exists an element $g\in \GL{n}{\Z[p^{-1}]}$ such that $g\G g^{-1} \leq \SL{n}{\Z[S^{-1}]}$.
\end{lemma}

\begin{proof}

Recall that all maximal compact subgroups of $\GL{n}{\Q_p}$ are conjugate and that $\GL{n}{\Z_p} \leq \GL{n}{\Q_p}$ is one such subgroup. The fact that it is both compact and open means that $\B_v := \GL{n}{\Q_p}/\GL{n}{\Z_p}$ is discrete. (The notation $\B_v$ is intended to remind the reader familiar with the Bruhat-Tits building for $\GL{n}{\Q_p}$ that $\B_v$ is the vertex set of the building, though we will not need to make use of that here.)

Also recall that the subgroup $\GL{n}{\Z[p^{-1}]} \leq \GL{n}{\Q_p}$ is dense, and since $\B_v$ is discrete, it follows that $\B_v = \GL{n}{\Z[p^{-1}]}/\GL{n}{\Z}$. (Observe that $\GL{n}{\Z} = \GL{n}{\Z[p^{-1}]} \cap \GL{n}{\Z_p}$.)

Now since the maximal compact subgroups of $\GL{n}{\Q_p}$ are in one to one correspondence with $\B_v$, we see that if $K \leq \GL{n}{\Q_p}$ is a maximal compact subgroup, then there exists an element $g \in \GL{n}{\Z[p^{-1}]}$ such that $K = g^{-1}\GL{n}{\Z_p}g$.

So, if $\G \leq \SL{n}{\Z[S_0^{-1}]} \leq \GL{n}{\Q_p}$ is precompact then $\G \leq K$ for some maximal compact subgroup $K$ of $\GL{n}{\Q_p}$ and by the above argument, there exists an element $g \in \GL{n}{\Z[p^{-1}]}$ such that 
\begin{equation}
 g\G g^{-1} \leq \GL{n}{\Z_p} \cap \SL{n}{\Z[S_0^{-1}]} = \SL{n}{\Z[S^{-1}]}.  \nonumber 
\end{equation}
\end{proof}


Now note that conjugation, as in Lemma 5.1, amounts to a change of basis. It is clear that if $(\sdp{\G}{_{\f'} \Q_p^N}, \Q_p^N)$ has relative property (T) then so does $(\sdp{\G}{_{\f''} \Q_p^N}, \Q_p^N)$ where $\f''$ is a conjugate representation of $\f'$. So, by Lemma 5.1, after conjugating if necessary, we may assume that $\f'(\G) \leq \SL{N}{\Z[S^{-1}]}$ and that $(\sdp{\G}{_{\f'} \Q_p^N}, \Q_p^N)$ has relative property (T) for each $p \in S \cup \{\8\}$.

By Lemma 5.2 (below), we have that the following group pair has relative property (T):
\begin{equation}
\(\sdp{\G}{(\Prod{ p \in S \cup \{\8\}} \Q_p^N )}, \Prod{ p \in S \cup \{\8\}} \Q_p^N \)
\nonumber
\end{equation}.

Finally, recall that the diagonal embedding $\Z[S^{-1}]^N \subset \Prod{ p \in S \cup \{\8\}} \Q_p^N$ is a co-compact lattice embedding. And, since $\f'(\G) \leq \SL{N}{\Z[S^{-1}]}$ it follows that this lattice is preserved by $\G$. Therefore, $\sdp{\G}{\Z[S^{-1}]^N}$ is a lattice in $\sdp{\G}{(\Prod{p \in S \cup \{\8\}} \Q_p^N)}$. Since lattices of this type inherit relative property (T) \cite[Proposition 3.1]{Jolissaint} this means that $(\sdp{\G}{_{\f'}\Z[S^{-1}]^N}, \Z[S^{-1}]^N)$ has relative property (T). 
\end{proof}

In the above proof, we made use of the following handy lemma:

\begin{lemma}
Suppose that $\G$ is a group acting by automorphisms on two groups $V_1$ and $V_2$. If $(\sdp{\G}{V_1}, V_1)$ and $(\sdp{\G}{V_2}, V_2)$ both have relative property (T) then $(\sdp{\G}{(V_1 \times V_2)}, V_1 \times V_2)$ also have relative property (T).
\end{lemma}

This is a corollary to the following general fact. The reader may notice the similarity between it and an analogous well known result about groups with property (T) and exact sequences.

\begin{lemma}
Suppose that $0 \to A_0 \to A \to A_1 \to 0$ is an exact sequence and that $\G$ acts by automorphisms on $A$ and leaves $A_0$-invariant. If $(\sdp{\G}{A_0}, A_0)$ and $(\sdp{\G}{A_1}, A_1)$ have relative property (T) then so does $(\sdp{\G}{A}, A)$.
\end{lemma}

\begin{proof}
Let $\pi \: \sdp{\G}{A} \to \Un(\H)$ be a unitary representation with almost invariant unit vectors $\{v_n\} \subset \H$. Then the space of $A_0$-invariant vectors $\H_0$ is non-trivial. Let $P \: \H \to \H_0$ and $P^\bot \: \H \to \H_0^\bot$ be the corresponding orthogonal projections. Observe that, since $A_0 \norm \sdp{\G}{A}$, the subspaces $\H_0$ and $\H_0^\bot$ are $\sdp{\G}{A}$-invariant and the corresponding projections commute with $\pi(\sdp{\G}{A})$.

We claim that for $n$ sufficiently large $\|P(v_n)\|^2 \geq 1/2$. Otherwise, there is a subsequence $n_j$ such that $\|P^\bot(v_{n_j})\|^2 = 1 - \|P(v_{n_j})\|^2 > 1/2$. Then
\begin{eqnarray}
\| \pi(\g)P^\bot(v_{n_j}) - P^\bot(v_{n_j}) \|^2 & = &\| P^\bot(\pi(\g)v_{n_j} - v_{n_j}) \|^2 \nonumber\\
\leq \| \pi(\g)v_{n_j} - v_{n_j} \|^2 & < & 2 \| \pi(\g)v_{n_j} - v_{n_j} \|^2 \cdot \|P^\bot(v_{n_j})\|^2  \nonumber
\end{eqnarray}
This of course means that if $v_{n_j}$ is $(K, \e)$ invariant then $P^\bot(v_{n_j})$ is $(K, \sqrt{2}\e)$-invariant. So, $\{P^\bot(v_{n_j})\} \in \H_0^\bot$ is a sequence of almost-invariant vectors, which is of course a contradiction: Indeed, $\H_0^\bot$ does not contain $A_0$-invariant vectors, so it can not contain $\sdp{\G}{A_0}$-almost invariant vectors. 

Therefore, for $n$ sufficiently large, $\|P(v_n)\|^2 \geq 1/2$. The same argument above shows that the restricted homomorphism $\pi_0 \: \sdp{\G}{A} \to \Un(\H_0)$ has almost invariant vectors $\{P(v_n)\}$. And since this homomorphism factors through $\sdp{\G}{A_1}$ we obtain the existence of a nonzero $A_1$-invariant vector. 
\end{proof}


\textbf{Remark:} Lemma 5.1 can be obtained in two other ways. One is a similar argument appealing to the CAT(0) structure of the Bruhat-Tits building for $\GL{n}{\Q_p}$ via a center of mass construction. Another is to observe that two maximal compact-open subgroups of $\GL{n}{\Q_p}$ are commensurable in the sense that their common intersection is a finite index subgroup in each. So, we may assume the result after passing to a finite index subgroup of $\G$.

\section{Algebro-Geometric Specialization}

In order to prove Theorem 1, in the direction of $(1) \implies (2)$, we need two basic ingredients. The first is to use the hypothesis (i.e. finite generation and the existence of a linear representation whose image has a non-amenable $\Re$-Zariski closure) in order to cook up a rational (or algebraic) representation to which we can apply Theorem 3, which is of course the second ingredient. This section is devoted to finding such a specialization, which is provided by the following:

\begin{prop}
Let $\G$ be a finitely generated group. If there exists a linear representation $\f \: \G\to \SL{n}{\Re}$ such that the Zariski closure $\~{\f(\G)}^Z(\Re)$ is non-amenable then there exists a representation $\psi \: \G \to \SL{m}{\Q}$ (possibly in a higher dimension) so that the Zarski closure $\~{\psi(\G)}^Z(\Re)$ is semisimple and not compact.
\end{prop}

Recall that a semisimple $\Re$-algebraic group is amenable if and only if it is compact. This follows from Whitney's theorem \cite[Theorem 3]{Whitney} (which says that a $\Re$-algebraic group has finitely many components as a $\Re$-Lie group) and from \cite[Corollary 4.1.9]{Zimmer} which states that a connected semisimple $\Re$-Lie group is amenable if and only if it is compact. So, the proposition guarantees that we may find, from an arbitrary $\Re$-representation, a $\Q$-representation which preserves the property of having non-amenable $\Re$-Zariski closure. The techniques used in the proof of this proposition are standard: the restriction of scalars functor and specializations of purely transcendental rings over $\Q$. 

However, we will also need a criterion which can distinguish when the image of a representation has non-amenable $\Re$-Zariski closure. This is provided by the following:

\begin{prop}
Let $\G$ be a finitely generated group. Then, there exists a normal finite index subgroup $\G_n \norm \G$ so that for any homomorphism $\f \:  \G \to \GL{n}{\Re}$ the following are equivalent:

\break
\begin{enumerate}
  \item The $\Re$-Zariski closure $\~{\f(\G)}^Z(\Re)$ is amenable.
  \item The traces of the commutator subgroup $\f([\G_n, \G_n])$ are uniformly bounded; that is
\begin{equation}
 |\tr(\f([\G_n, \G_n]))| \leq n. \nonumber
\end{equation}
\end{enumerate}
\end{prop}

\textbf{Remark:} It is a fact (see Subsection 6.3, Lemma 6.6), that if a subgroup of $\GL{n}{\Re}$ has bounded traces, then it's $\Re$-Zariski closure is amenable (actually it is a compact extension of a unipotent group). Therefore, in the direction of (2) implies (1), there is nothing special about $[\G_n, \G_n]$. Namely, any co-amenable normal subgroup of $\G$ would do. The more subtle direction is that of (1) implies (2). It is in this direction that we must work to find a suitable $\G_n$. Under the added assumption that $\~{\f(\G)}^Z(\Re)$ is Zariski-connected the result follows from classical structure theory of Zariski-connected $\Re$-algebraic groups with $\G_n = \G$. 

However, we must address the fact that the image of a general representation $\f \: \G \to \GL{n}{\Re}$, need not have Zariski-connected Zariski-closure. It turns out that for an arbitrary (reductive) $\Re$-algebraic group, there is a finite index subgroup (with uniformly bounded index) which ``behaves as if'' it were connected (see Subsection 6.2, Lemma 6.3). Namely, it has most of the nice structure properties of Zariski-connected groups (see Subsection 6.1, Lemma 6.2). It turns out that the uniform bound on the index of this subgroup, together with its ``pseudo-connectedness'' properties are exactly what we need to find a suitable $\G_n$ which is done in Subsection 6.4. We then prove Proposition 6.2 in Subsection 6.5 and Proposition 6.1 in Subsection 6.6.

\subsection{Some Algebraic Facts}

Throughout this section, we will be dealing exclusively with $\Re$-Zariski closures. As such we will write $G$ instead of $G(\Re)$, when speaking of $\Re$-Zariski closed groups, and we will just say Zariski-closed or algebraic. Also, when we say connected, we mean Zariski-connected. We now develop the necessary lemmas to prove Proposition 6.2. 

\begin{definition}
An algebraic group $G$ is said to be reductive if any closed unipotent normal subgroup is trivial. 
\end{definition}

Observe that it is common to require in the definition of a reductive group that either $G$ be Zariski-connected or that any closed \emph{connected} normal unipotent subgroup of $G$ be trivial. 

However, in characteristic zero, the two notions are the same since algebraic unipotent groups are always Zariski-connected. This follows by 
\begin{itemize}
  \item Chevalley's Theorem: \cite[Theorem 11.2]{Humphreys} If $H \leq G$ are two algebraic groups, then there exists a rational representation $G \to \GL{N}{\Re}$ and a vector $v \in \Re^n$ such that $H = \stab_G(\Re\cdot v)$.
  \item The image of a unipotent element under a rational homomorphism is unipotent.
  \item Unipotent elements have infinite order in characteristic zero.
\end{itemize}  

\break
To be complete, we also give the following definition:

\begin{definition}
An algebraic group $G$ is said to be semisimple if any closed solvable normal subgroup is finite. 
\end{definition}

And now onto the lemmas; the first of which shows that we may restrict our attention to reductive groups, since doing so does not affect the hypotheses and conclusions of Proposition 6.2. 


\begin{lemma}
Suppose that $L \leq \GL{n}{\C}$ is a $\C$-closed group and $U \norm L$ is the maximal unipotent normal subgroup. There is a representation $\pi \: L \to \GL{n}{\C}$ such that $\ker(\pi) = U$ and $\tr(g) = \tr(\pi(g))$ for every $g \in L$.
\end{lemma}

\begin{proof}
Choose a Jordan-H\"older series for $\C^n$ as an $L$-module:
\begin{equation}
0 = V_0 \subset V_1\subset \cdots \subset V_k = \C^n
\nonumber
\end{equation} 
so that the corresponding representation $\rho_i \: L \to \GL{}{V_i/V_{i-1}}$ is irreducible. Then for each $i$, the image $\rho_i(U)$ is again unipotent, and by the Lie-Kolchin Theorem there is a vector $v_i \in V_i/V_{i-1}$ which is fixed by $\rho_i(U)$. But since $\rho_i(U)\norm \rho_i(L)$ and $\rho_i(L)$ acts irreducibly on $V_i/V_{i-1}$ it follows that $U \leq \ker(\rho_i)$. 

Choosing a basis, for $\C^n$ which respects this Jordan-H\"older series, we see that 
\begin{equation}
L \leq %
\(
\begin{array}{cccc}
    \rho_1(L)  & * & \cdots & *  \\
      0 & \rho_2(L) & \cdots & * \\
      \vdots & \vdots & \ddots & \vdots \\
      0 & 0 & \cdots & \rho_k(L) \\      
\end{array}
\)
\nonumber
\end{equation}
and
\begin{equation}
U \leq %
\(
\begin{array}{cccc}
    I_{n_1}  & * & \cdots & *  \\
      0 & I_{n_2} & \cdots & * \\
      \vdots & \vdots & \ddots & \vdots \\
      0 & 0 & \cdots & I_{n_k} \\      
\end{array}
\),
\nonumber
\end{equation}
where $n_i$ is the dimension of $V_i/V_{i-1}$ and $I_{n_i}$ is the $n_i \times n_i$ identity matrix. 

Let $V = \underset{i = 1}{\overset{k}{\bigoplus}} V_i/V_{i-i}$ be the corresponding $n$-dimensional vector space. Consider the homomorphism $ \pi \: L \to \underset{i = 1}{\overset{k}{\bigoplus}} \rho_i (L) \leq \GL{}{V}$. By construction, $\tr(\pi(g)) = \Sum{i = 1}{k} \tr(\rho_i(g)) = \tr(g)$ for each $g \in L$. Furthermore, it is clear that $\ker(\pi)$ is unipotent, and contains $U$. Since $U$ is maximal it follows that $\ker(\pi) = U$.
\end{proof}


The following is a corollary to the proof above:

\begin{cor}
Let $G$ be a $\Re$-algebraic reductive group. Then every $\C$-representation of $G$ is the direct sum of $G$-irreducible sub-representations. 
\end{cor}

This next lemma is classical. These are exactly the ``nice'' properties of connected (and reductive) groups that were alluded to above.

\begin{lemma}
Let $G_0$ be a connected reductive group. Then the following hold: 
\begin{enumerate}
  \item The radical $R(G_0) = Z(G_0)^o$, where $Z(G_0)^o$ is the identity component of the center of $G_0$.
  \item The intersection $[G_0, G_0] \cap Z(G_0)$ is finite.
  \item $G_0 = [G_0, G_0] \cdot Z(G_0)$.
  \item The commutator subgroup $[G_0, G_0]$ is semisimple.
\end{enumerate}
\end{lemma}

\begin{proof}
For assertions (1) and (2) we cite \cite[Lemma 19.5]{Humphreys}. 

Assertion (3) follows from (2) by noting that $G_0/[G_0, G_0] \cdot Z(G_0)$ is a connected Abelian semisimple group, and therefore trivial. 

Assertion (4) follows from (3) and (2): Let $R \norm [G_0, G_0]$ be a closed solvable normal subgroup. Since $G_0$ is reductive, $G_0/R(G_0)$ is semisimple. Then, $R/R\cap R(G_0)$ is closed and solvable and hence finite. Since $[G_0, G_0] \cap R(G_0)$ is finite, it follows that $R$ is finite.
\end{proof}

This next lemma yields the want-to-be connected group that was alluded to above.

\begin{lemma}
Let $G_0$ be a connected reductive group of finite index in $G \leq \GL{n}{\Re}$. Then there exists a subgroup $G_1 \norm G$ such that
\begin{enumerate}
  \item $G_0 \norm G_1$.
  \item The index $[G:G_1] \leq n!$.
  \item The commutator subgroup $[G_1, G_1]$ contains $[G_0, G_0]$ as a finite index normal subgroup. (And hence $[G_1, G_1]$ is semisimple.)
\end{enumerate}
\end{lemma}

We first prove the following special case:

\begin{lemma}
Let $G_0$ be a connected reductive group of finite index in $G$. Suppose that $G \leq \GL{n}{\C}$ is an irreducible representation. Then, there exists a subgroup $G_1 \norm G$ with the following properties:

 \begin{enumerate}
 
  \item $G_0 \norm G_1$
  \item The index $[G : G_1] \leq n!$.
  \item If $Z(G_0)$ and $Z(G_1)$ are the centers of $G_0$ and $G_1$ respectively then $Z(G_0) \leq Z(G_1)$.
\end{enumerate}
\end{lemma}

\begin{proof}

Since $G_0$ is reductive, the representation on $\C^n$ decomposes as a direct sum of irreducible sub-representations. Let $V \subset \C^n$ be one such.

Now, since $G_0 \norm G$ it follows that for each $g \in G$ the subspace $gV$ is also an irreducible $G_0$-sub-representation. Hence if $gV \cap V \neq \{0\}$ then $gV = V$. 

\begin{clam}
There exists $\{g_1, \dots, g_l\} \subset G$ such that $\C^n = \Oplus{j = 1}{l} g_jV$.
\end{clam}

\begin{proof}
Let $g_1 = 1$. Then either $V= \C^n$, or there exists a $g_2 \in G$ such that $V \cap g_2 V = \{0\}$. In this latter case we have that $V\oplus g_2 V \subset \C^n$. 

Inductively, suppose that we have found $\{g_1, \dots, g_k\} \subset G$ such that the corresponding $g_jV$ are linearly independent. Namely so that $\Oplus{j=1}{k} g_jV \subset \C^n$ is a direct sum of $G_0$-irreducible sub-representations.

Observe that $\Oplus{j=1}{k} g_jV$ is $G_0$-invariant. And since the $G$-translates of $V$ are $G_0$-irreducible sub-representations we get the following dichotomy:

\begin{enumerate}
  \item There exists a $g_{k+1} \in G$ such that $g_{k+1}V \cap \Oplus{j=1}{k} g_jV = 0$, or
  \item $gV \subset \Oplus{j=1}{k} g_jV$ for each $g \in G$.
\end{enumerate}

In case (1) we may conclude that $\Oplus{j=1}{k+1} g_jV \subset \C^n$ is a direct sum of $G_0$-irreducible sub-representations.

In case (2) we must have that $gg_iV \subset \Oplus{j=1}{k} g_jV$ for each $i = 1, \dots, k$, and $g \in G$. This means that $\Oplus{j=1}{k} g_jV$ is $G$-invariant, and hence $\Oplus{j=1}{k} g_jV = \C^n$.

\end{proof}
 
This induces a homomorphism $\s \: G \to \textrm{Sym}(l)$ where the $\textrm{Sym}(l)$ denotes the symmetric group on $l$-symbols. Let $G_1 = \Cap{j = 1}{l} \stab_G(g_jV)$. Then clearly, $G_1 = \ker(\s)$, so that $G_1$ satisfies properties (1) and (2) as promised above. (Note that $l \leq n$.) 

Furthermore, all of the $G_0$-irreducible subspaces are $G_1$-invariant and hence these are also $G_1$-irreducible subspaces. 
By Schur's Lemma, the centers of $G_0$ and $G_1$ are block-scalar matrices of the same type, and therefore, $G_1$ also satisfies property (3) as it was promised to do. 
\end{proof}

In order to pass from Lemma 6.4 to Lemma 6.3 we will need the following:

\begin{lemma}
If $G_0 \norm G$ is a finite index subgroup then $[G, G_0]$ is a normal finite index subgroup of $[G, G]$. 
\end{lemma}

\begin{proof}
Since $G_0 \norm G$ it follows that $[G, G_0] \norm G$ (and in particular $[G,G_0] \norm [G,G]$). Hence, to show that the index of $[G, G_0]$ in $[G,G]$ is finite, it is sufficient to show that if $[G, G_0] =1$ then $[G,G]$ is finite. (Just take the quotient of $G$ by $[G, G_0]$ if necessary, and use the general fact that for any homomorphism $h \: G \to H$ and any subgroups $A , B \leq G$ the following equality holds: $h([A,B]) = [h(A),h(B)]$.)

If $[G, G_0] =1$, it follows that $G$ centralizes $G_0$. That is, $G_0 \leq Z(G)$. Then,
\begin{equation}
[G:Z(G)] \leq [G:G_0] <\8
\nonumber
\end{equation} 

This implies that $[G,G] <\8$ (see \cite[Lemma 17.1.A]{Humphreys}).
\end{proof}


\subsection{Proof of Lemma 6.3}

\begin{proof}
The assumptions are that $G_0 \norm G \leq \GL{n}{\Re}$ where $G_0$ is connected, reductive and of finite index in $G$. This means that $G$ is also reductive and so by Corollary 6.1 we have that the representation of $G$ on $\C^n = \Oplus{i = 1}{k} \C^{n_i}$ is the direct sum of $G$-irreducible subrepresentations. By considering each irreducible piece and applying Lemma 6.4, we see that there exists a subgroup $G_1 \norm G$ of index at most $\overset{k}{\Prod{i = 1}} (n_i)! \leq n!$ such that $G_0 \norm G_1$ and $Z(G_0) \leq Z(G_1)$.

We claim that $[G_1, G_0] = [G_0, G_0]$:

Let $ x \in G_1, y \in [G_0, G_0]$, and $z \in Z(G_0) \leq Z(G_1)$. Recall that $[G_0, G_0] \norm G_1$ so that

\begin{equation}
[x,yz] = [x,y] = (xyx^{-1})y^{-1} \in [G_0, G_0]
\nonumber
\end{equation}

Since $G_0 = [G_0, G_0] \cdot Z(G_0)$ it follows that $[G_1, G_0]$ is generated by elements in $[G_0, G_0]$ and therefore, $[G_1, G_0] \leq [G_0, G_0]$. On the other hand, $[G_0, G_0] \leq [G_1, G_0]$ so $[G_1, G_0] = [G_0, G_0]$. 

Now, since $G_0$ has finite index in $G_1$ by Lemma 6.5, we see that $[G_0, G_0] = [G_1, G_0]$ is a finite index normal subgroup of $[G_1, G_1]$ and we are done.
\end{proof}


\subsection{The Trace Connection}

So, far, we have addressed only the structure of the algebraic groups in question, and have ignored the role of the trace. We now discuss how the trace ties in to the picture.

Recall that if $\G \leq \GL{n}{\Re}$ is a precompact group then all of its eigen values have norm 1 and hence its traces are uniformly bounded by $n$. Also recall that the Zariski closure of a precompact group is compact and therefore amenable. The following shows that the converse also holds. Namely:

\begin{lemma}
Let $\G \leq \GL{n}{\Re}$ be a group. If the set of traces $\tr(\G) := \set{\tr(\g)}{\g \in \G}$ is bounded then the Zariski-closure $\~{\G}^Z(\Re)$ is amenable.
\end{lemma}

\subsubsection{Some Useful Facts}

We will need the following:

\begin{fact} \cite[Corollary 1.3(c)]{Bass}
Let $\G \leq \GL{}{V}$ be a group acting irreducibly on the complex vector space $V$. If the traces of $\G$ are bounded then $\G$ is precompact (in the $\C$-topology).
\end{fact}

\begin{claim}
Let $\G \leq \GL{n}{\C}$ be a subgroup such that $B = \Sup{\g \in \G}|\tr(\g)| < \8$. Then all $\G$-eigen values have norm 1 and $B = n$.
\end{claim}

\begin{proof}

By contradiction suppose that there is some $\g \in \G$ with an eigen value of norm not equal to 1. Then upon passing to $\g^{-1}$ if necessary, we may assume that $\g$ has an eigen value of norm strictly greater than 1. 

Order the eigen values so that
 \begin{equation}
|\lambda_1| = \cdots =|\lambda_m| > |\lambda_{m+1}| \geq \cdots |\lambda_n|. \nonumber
\end{equation} 

Since the traces of $\G$ are bounded we get that for each $k \in \N$

\begin{equation}
\l|\tr(\g^k)\r| = \l|\Sum{j=1}{n} \lambda_j^k \r| \leq B.
\nonumber
\end{equation}

\noindent
The triangle inequality gives us that

\begin{equation}
\l|\Sum{j=1}{m} \frac{\lambda_j^k}{\lambda_1^k} \r| \leq \frac{ B}{|\lambda_1^k|} +\Sum{j=m+1}{n}\l| \frac{\lambda_j^k}{\lambda_1^k} \r| \to 0.
\nonumber
\end{equation}

\noindent
By Claim 6.2 (see below) we get that $\Sum{j=1}{m} \frac{\lambda_j^k}{\lambda_1^k} = m$ and so 

\begin{equation}
1 \leq m = \l|\Sum{j=1}{m} \frac{\lambda_j^k}{\lambda_1^k} \r|  \to 0
\nonumber
\end{equation}
a contradiction.

Therefore, all eigen values of $\G$ have norm 1 and the supremum $B = \Sup {\g \in \G}|\tr(\g)|$ is attained at the identity.

\end{proof}

\begin{claim} If $\Sum{j =1}{n} e^{ik\theta_j}$ converges as $k \to \8$ then $\Sum{j =1}{n} e^{ik\theta_j} \equiv n$ for all $k$.

\end{claim}

\begin{proof}
Consider the action of $\Z$ by the rotation on the $n$-torus $\T^n$ corresponding to $(e^{i\theta_1}, \dots, e^{i\theta_n})$. Let us definte the closed subgroup
\begin{equation}
S = \~{\<(e^{ik\theta_1}, \dots, e^{ik\theta_n}) \: k \in \Z\>}.
\nonumber
\end{equation} 

Now, if $S$ is discrete then it is finite, which means that the identity $(1, \dots, 1) \in \T^n$ is a periodic point. On the other hand, if $S$ is not discrete then the identity is an accumulation point of the sequence $\{(e^{ik\theta_1}, \dots, e^{ik\theta_n}) \: k \in \Z\}$. Either way, there is a subsequence $k_l \to +\8$ such that 
\begin{equation}
\Lim{l \to \8}(e^{i{k_l}\theta_1}, \dots, e^{i{k_l}\theta_n}) = (1, \dots, 1).
\nonumber
\end{equation}

This shows that if $\Sum{j =1}{n} e^{ik\theta_j}$ converges then 
\begin{equation}
\Lim{k \to \8} \Sum{j =1}{n} e^{ik\theta_j} = \Lim{l \to \8} \Sum{j =1}{n} e^{i{k_l}\theta_j} = n.
\nonumber
\end{equation}

Since the sequence is convergent, any subsequence converges to the same limit. Therefore the same argument shows that if $(e^{i\psi_1}, \dots, e^{i\psi_n}) \in S$ then  

\begin{equation}
\Sum{j = 1}{n}e^{i\psi_j} = n.
\nonumber
\end{equation}

In particular, this holds for $(e^{i\psi_1}, \dots, e^{i\psi_n}) = (e^{i\theta_1}, \dots, e^{i\theta_n})$. Since 1 is an extreme point of the unit disk we conclude that $e^{i\theta_j} = 1$ for each $j =1, \dots, n$.

\end{proof}

\subsubsection{The Proof of Lemma 6.6}

\begin{proof}
By Lemma 6.1 we may assume that $G = \~{\G}^Z(\Re) \leq \GL{n}{\Re}$ is reductive. Using Corollary 6.1, we decompose $\C^n = \Oplus{i \in I}{}V_i$ into a direct sum of $G$-irreducible sub-representations. Since $\G$ is Zariski-dense in $G$, this of course means that each $V_i$ is also a $\G$-irreducible sub-representation. We aim to show that $G$ is compact. To this end, it is sufficient to show that $\G$ is pre-compact in $\GL{n}{V_i}$  for each $i \in I$ since the homomorphism $G \to \Prod{i \in I} \GL{n}{V_i}$ is rational and injective.

By Claim 6.1, $\G$ has bounded traces in each $\GL{n}{V_i}$.
And by Fact 6.1, $\G$ is precompact in each $\GL{n}{V_i}$ since it acts irreducibly on $V_i$. 
\end{proof}

\subsection{Choosing $\G_n$ for Proposition 6.2}

Recall that condition (2) of Lemma 6.3 guarantees a uniform bound on the index of the groups in question. We now show how we will make use of that fact to find our $\G_n$:

\begin{lemma}
Let $\G$ be a finitely generated group and let
\begin{equation}
H_N \: \{\pi \: \G \to F| F \text{ is a group of order at most }N\}.
\nonumber
\end{equation}
 Then $\G(N) = \Cap{\pi \in H_N}{} \ker(\pi)$ is a finite index (normal) subgroup of $\G$. 
\end{lemma}

\begin{proof}
This is a straightforward consequence of two facts:

\begin{enumerate}
  \item[Fact 1:] There are finitely many groups of order at most $N$. 
  \item[Fact 2:] There are finitely many homomorphisms from a finitely generated group to a fixed finite group.  
\end{enumerate}
\end{proof}


\subsection{The proof of Proposition 5}

\begin{proof}
Let $\G_n = \G(n!)$ as in Lemma 6.7.  Then, $\G_n$ is a finite index normal subgroup of $\G$. Let $\f \: \G \to \GL{n}{\Re}$ be any homomorphism. 

$(2) \implies (1)$: 
If the set of traces $\tr(\f([\G_n, \G_n]))$ is uniformly bounded, then by Lemma 6.6 the Zariski closure $\~{\f([\G_n,\G_n])}^Z(\Re)$ is amenable. Therefore, $\~{\f(\G)}^Z(\Re)$ is amenable as it is a virtually Abelian extension of $\~{\f([\G_n,\G_n])}^Z(\Re)$.

$(1) \implies (2)$: 
Suppose that $G := \~{\f(\G)}^Z(\Re)$ is amenable. As was mentioned several times, by Lemma 6.1 it does no harm to assume that $G$ is reductive. Let $G_0$ be the Zariski connected component of 1 and let $G_1$ be as in Lemma 6.3. Then, $[G_0, G_0]$ being Zariski-connected, semisimple, and amenable, it is compact. Since $[G_1, G_1]$ contains $[G_0, G_0]$ as a finite index subgroup, it follows that $[G_1, G_1]$ is compact. 

Thus, if $\f(\G_n) \leq G_1$ then we are done. But, this follows by construction: Recall that $\G_n \leq \ker(\pi)$ for every homomorphism $\pi \: \G \to F$ where $F$ is a finite group of order at most $n!$. Since $G_1 \norm G$ and the index $[G:G_1] \leq n!$ we must have that $\f(\G_n) \leq G_1$. 

\end{proof}

\subsection{Finally: The Proof of Proposition 6.1}

\begin{proof}

To conserve notation, we assume that $\G \leq \SL{n}{\Re}$. Let $K$ be the field generated by the entries of some finite generating set for $\G$ so that $\G \leq \SL{n}{K}$. Then, since $K$ is finitely 

generated, it is a finite and hence separable extension of $\Q(t_1, \dots, t_s) \subset \Re$, where $t_1, \dots, t_s \in K$ are algebraically independent transcendentals. So, after applying the restriction of scalars if necessary, we may assume that $\G \leq \SL{n}{\Q(t_1, \dots, t_s)}$. (We note that property (3) of the restriction of scalars, guarantees that the hypothesis is preserved.) 

The proof is by induction on the transcendence degree of $\Q(t_1, \dots, t_s)/\Q$. 

\textbf{Base Case:} Suppose $s=0$. 

Let $G = \~{\G}^Z$ be the Zariski-closure. Since $\G \leq \SL{n}{\Q}$ it follows that $G$ and its radical $R(G)$ are defined over $\Q$. Fixing a representation of $G/R(G) (\Q) \leq \SL{n}{\Q}$ we have the desired result.

\textbf{Induction Hypothesis:} Assume it is true for $s-1$. 

Since $\G$ is finitely generated, it follows that there exist irreducible polynomials $\d_1, \dots, \d_l \in \Q[t_1, \dots, t_s]$ such that if we set $\R = \Q[t_1, \dots, t_s, \d_1^{-1}, \dots, \d_l^{-1}]$ then $\G \leq \SL{n}{\R}$. 

Observe that by Proposition 6.2, $[\G_n, \G_n] \leq \SL{n}{\Re}$ has unbounded traces since $\~{\G}^Z(\Re)$ is non-amenable. So, up to a relabeling of the transcendentals there are two cases to consider:

\begin{enumerate}
  \item [Case 1:] The unbounded traces of $[\G_n,\G_n]$ are independent of $t_s$, that is 
  \begin{equation}
\set{\tr(\g)}{\g \in [\G_n,\G_n] \text{ and } |\tr(\g)| \geq n+2} \cap \R \subset \Q(t_1, \dots, t_{s-1}).
\nonumber
\end{equation}
  \item [Case 2:] There is an element in $[\G_n,\G_n]$ with large trace which is non-constant as a rational function in $t_s$. Namely, there is a $\g \in [\G_n,\G_n]$ such that $|\tr(\g)| \geq n+2$ and $\tr(\g) \in \R \bs \Q(t_1, \dots, t_{s-1})$.
\end{enumerate}

We now need to say how we will specialize the transcendental $t_s$. First consider the denominators $\d_i$ as polynomials in $t_s$. Since there are finitely many, the bad set 
\begin{equation}
B = \set{ a \in \Re}{\d_i(t_1, \dots, t_{s-1}, a) = 0 \text{ for some } i = 1, \dots l}
\nonumber
\end{equation} 
is finite. (Recall that the kernel of a ring homomorphism cannot contain any invertible elements.) Now, we choose the specialization in each case:

\break
\begin{enumerate}
  \item [Case 1:] Choose $a \in \Q \bs B$.
  \item [Case 2:] Let $\g \in [\G_n,\G_n]$ be such that 
 \begin{equation}
r(t_s)=\tr(\g)
\nonumber
\end{equation}
is a nonconstant rational function in $t_s$ and such that $|r(t_s)|\geq n+2$. Then, $r(x)$ is a continuous function in some neighborhood of $t_s \in \Re \bs B$ and so there is an $a \in \Q \bs B$ such that $|r(a)| \geq n+1$.
\end{enumerate}

Now, fix an embedding $\Q(t_1, \dots, t_{s-1}) \subset \Re$ and let
\begin{equation}
 \psi \: \SL{n}{\Q(t_1, \dots, t_s)} \to \SL{n}{\Q(t_1, \dots, t_{s-1})} \nonumber
\end{equation}
be the homomorphism induced from the ring homomorphism $t_s \mapsto a$. Observe that this is well defined since we are dealing with unimodular matrices. 

To apply the induction hypothesis, we must show that the Zariski-closure $\~{\psi(\G)}^Z(\Re)$ is again non-amenable. This is immediate by Proposition 6.2 since by construction, there is a $\g \in [\G_n,\G_n]$ such that $|\tr(\psi(\g))| \geq n+1$. Since the traces of a subgroup of $\SL{n}{\Re}$ are either uniformly bounded by $n$ or unbounded, we see that $\psi([\G_n,\G_n])$ has unbounded traces and the proposition is proved. 
\end{proof}


\section{Proof of Theorem 1 in the direction $(1) \implies (2)$}

We instead prove the following:

\begin{theorem}
Suppose that $\G$ is a finitely generated group which admits a linear representation $\f \: \G \to \SL{n}{\Re}$ such that the $\Re$-Zariski closure $\~{\f(\G)}^Z(\Re)$ is non-amenable. Then there exists a finite set of primes $S \subset \Z$ and a homomorphism $\a \: \G \to \SL{N}{\Z[S^{-1}]}$ such that, if $A=\Z[S^{-1}]^N$ then $(\sdp{\G}{_\a A}, A)$ has relative property (T).
\end{theorem}


The proof is in two basic steps:
\begin{enumerate}
  \item [\textbf{Step A:}] Show that under the hypothesis of Theorem 7.1 there is a finite index subgroup $\G_0 \norm \G$ satisfying property (F$_\8$). 
  \item [\textbf{Step B:}] Show that if $\G_0 \norm \G$ is a finite index subgroup such that $(\G_0 \ltimes A, A)$ has relative property (T) then there is an action of $\G$ on $A^k$ (with $A$ as in Theorem 7.1 and $k = [\G:\G_0]$) such that $(\G \ltimes A^k, A^k)$ has relative property (T).
  
\end{enumerate}

It is clear that Steps A and B prove Theorem 7.1 by Theorem 3.

\subsubsection{Proof of Step A}

\begin{proof}

By Proposition 6.1 there exists a rational representation $\psi \: \G \to \SL{m}{\Q}$ such that $\~{\psi(\G)}^Z (\Re)$ is semisimple and not compact.

Let $\G_0 \norm \G$ be the normal subgroup of finite index such that $\~{\psi(\G_0)}^Z(\Re)$ is the Zariski-connected component of the identity of $\~{\psi(\G)}^Z(\Re)$. Then, $G(\Re):=\~{\psi(\G_0)}^Z(\Re)$ is again not compact semisimple. 

In order to be totally precise, we now turn our attention to the $\C$-Zariski closure $G(\C)$, which is of course defined over $\Q$. Furthermore, we fix an embedding $\~\Q \subset \C$. 

\textbf{Step A.1:} There is a $\~\Q$-homomorphism $\pi \: G(\C) \to \Prod{i \in I} H_i(\C)$ with finite central kernel, where each $H_i(\C)$ is a $\~\Q$-simple $\~\Q$-group.

Since $G(\C)$ is Zariski-connected and semisimple, this follows from \cite[Proposition 2]{Tits}. Let $\pi_i \: G(\C) \to H_i(\C)$ be the corresponding $\~\Q$-projection. 

\textbf{Step A.2:} Each $H_i$ is defined over $K_i$, a finite normal extension of $\Q$ and $\pi_i$ is a $K_i$-morphism.

By \cite[ Propositions 3.1.8 \& 3.1.10]{Zimmer}, this follows from the fact that $\pi_i \psi(\G_0) \leq H_i (\~\Q)$ is a Zariski-dense finitely generated subgroup. 

Now, for each $i$, fix a $K_i$-rational representation $H_i(\~\Q) \to \SL{n_i}{\~\Q}$ without fixed vectors and identify $H_i(\~\Q)$ with its image. By abuse of notation, we still take $\pi_i \: \G_0 \to H_i(K_i) \leq \SL{n_i}{\~\Q}$.

\textbf{Step A.2:} There is an $i_0$ such that $\pi_{i_0}$ realizes property (F$_\8$) for $\G_0$.

Observe that by construction, the $\~\Q$-Zariski-closure of $\pi_i(\G_0)$ is $H_i(\~\Q)$ and is therefor $\~\Q$-simple. For the same reason $\pi_i(\G_0)\leq \SL{n_i}{\~\Q}$ has no fixed vectors as this is a Zariski-closed condition. Thus in order for $\pi_i$ to realize property (F$_\8$) for $\G_0$ we need only show that the corresponding diagonal embedding into $R_{K_i/\Q}H_i(\Re)$ is not precompact. We now find an $i_0$ for which this holds.

Recall that the restriction of scalars satisfies several nice properties, which were enumerated in Section 2. We will refer to these by number below:

Let $i \in I$. Recall that by Property 1, the restriction of scalars $\R_{K_i/\Q}H_i(\C)$ is uniquely determined (up to $\Q$-isomorphism) by specifying a ``projection'' $P_i \: \R_{K_i/\Q}H_i(\C) \to H_i(\C)$, which we now fix. 

Since $G(\C)$ is a $\Q$-group and $\pi_i$ is a $K_i$-morphism, it follows (Property 2) that there is a unique $\Q$-morphism $\rho_i \: G(\C) \to \R_{K_i/\Q}H_i(\C)$ so that $\pi_i = P_i \circ \rho_i$.

This of course means that there is a $\Q$-morphism $\rho \: G(\C) \to \Prod{i \in I} \R_{K_i/\Q}H_i(\C)$ such that $\pi = P \circ \rho$ where $P \: \Prod{i \in I} \R_{K_i/\Q}H_i(\C) \to \Prod{i \in I} H_i(\C)$ is the obvious projection. Furthermore, the kernel of $\rho$ is finite since $\ker(\rho) \leq \ker(\pi)$. So, $\rho$ is virtually an isomorphism onto its image. 

Now, since $\rho$ is a $\Q$-morphism with finite kernel, it follows that $\rho( G(\Re)) \leq \Prod{i \in I} \R_{K_i/\Q}H_i(\Re)$ is semisimple and not compact. This means that for some $i_0 \in I$ the corresponding homomorphism $\rho_{i_0} \: \psi(\G_0) \to \R_{K_{i_0}/\Q}H_{i_0}(\Q)$ has non-precompact image in $\R_{K_{i_0}/\Q}H_{i_0}(\Re)$. 
\end{proof}

\subsubsection{Proof of Step B}

\begin{proof}

Let $\a \: \G_0 \to \SL{N}{\Z[S^{-1}]}$ such that, setting $A = \Z[S^{-1}]^N$, we have that  $(\sdp{\G_0}{_{\a} A}, A)$ has relative property (T). 
  
Also, let $k = [\G : \G_0]$. We now construct a homomorphism  $\a' \: \G \to \SL{kN}{\Z[S^{-1}]}$ such that $( \sdp{\G}{A^k}, A^k )$ has relative property (T). 

Set $F = \G/\G_0$ and choose a section $s \: F \to \G$ such that if $[\cdot] \: \G \to F$ is the natural projection then for any $f \in F$, $[s(f)] = f$. Let $c \: \G \times F \to \G_0$ be the corresponding cocycle. That is, $c(\g, f) = s(\g f)^{-1} \g s(f)$. 

Define the action of $\G$ on $\Oplus{f \in F}{} A$ as follows: 
\begin{equation}
\g(a_f)_{f \in F} = (c(\g, f) \cdot a_f)_{\g f \in F}
\nonumber
\end{equation}

The fact that $c$ is a cocycle ensures that this is a well defined action, and it is clearly by automorphisms since $\G_0$ acts by automorphisms. Therefore, we may form the semidirect product $\sdp{\G}{\Oplus{f \in F}{} A}$. 

To show that $(\sdp{\G}{\Oplus{f \in F}{} A}, \Oplus{f \in F}{} A)$ has relative property (T) it is sufficient to show that \linebreak
$(\sdp{\G_0}{\Oplus{f \in F}{} A}, \Oplus{f \in F}{} A)$ has relative property (T). Indeed, any unitary representation of $\sdp{\G}{\Oplus{f \in F}{} A}$ is a (continuous) unitary representation of $\sdp{\G_0}{\Oplus{f \in F}{} A}$.

Now, observe that since $\G_0 \norm \G$ the corresponding $\G_0$ action on $\Oplus{f \in F}{} A$ is given by: 
\begin{equation}
\g_0(a_f)_{f \in F} = (s(f)^{-1}\g_0 s(f) \cdot a_f)_{f \in F}.
\nonumber
\end{equation}
Namely, $\G_0$ preserves the $f_0$-component $A_{f_0} \leq \Oplus{f \in F}{} A$ for each $f_0 \in F$.

Let $\sdp{\G_0}{_{s(f_0)}A} \leq \sdp{\G_0}{\Oplus{f \in F}{} A}$ be the subgroup corresponding to $f_0 \in F$. It follows from Lemma 5 that if $(\sdp{\G_0}{_{s(f_0)}A}, A)$ has relative property (T) for each $f_0 \in F$ then $(\sdp{\G_0}{\Oplus{f \in F}{} A}, \Oplus{f \in F}{} A)$ has relative property (T). 

And this is indeed the case since twisting the $\G_0$-action by $s(f_0)$ amounts to precomposing the $\G_0$-action on $A$ by an automorphism of $\G_0$. And, the conclusion of Burger's Criterion, and hence the proof of Theorems 2 and 3, remains valid under this twist.
\end{proof}



 \section{Theorem 1 in the direction of $(2) \implies (1)$}

Recall that there is a natural embedding $\GL{n}{\Re} \leq \SL{n+1}{\Re}$ induced by 
\begin{equation}
g \mapsto \mathrm{diag}(g, 1/\mathrm{det}(g)). \nonumber
\end{equation}
Hence, $\SL{n}{\Re} \leq \GL{n}{\Re} \leq \SL{n+1}{\Re}$. This means that there is a homomorphism $\f \: \G \to \GL{n}{\Re}$ such that $\~{\f(\G)}^Z(\Re)$ is non-amenable if and only if there is a homomorphism $\f' \: \G \to \SL{n'}{\Re}$ such that $\~{\f'(\G)}^Z(\Re)$ is non-amenable. This shows that Theorem 1 is equivalent to the following:
\break
\begin{theorem1'}
Let $\G$ be a finitely generated group. The following are equivalent:
\begin{enumerate}
  \item There exists a homomorphism $\f : \G \to \GL{n}{\Re}$ such that the $\Re$-Zariski-closure $\~{\f(\G)}^Z(\Re)$ is non-amenable.
  \item There exists an Abelian group $A$ of nonzero finite $\Q$-rank and a homomorphism $\f' \: \G \to \Aut(A)$ such that the corresponding group pair $(\sdp{\G}{_{\f'}A}, A)$ has relative property (T).
\end{enumerate}
\end{theorem1'}

So, in this section, we will show Theorem 1' in the direction of $(2) \implies (1)$. To do this we will make use of the following simple facts:

\begin{lemma}
Suppose that $\G$ is finitely generated and $A$ is countable. If $(\sdp{\G}{A}, A)$ has relative property (T) then $\sdp{\G}{A}$ is finitely generated.
\end{lemma}

The proof of this is exactly as one would show that a countable group with property (T) is finitely generated. See \cite[Theorem 7.1.5]{Zimmer}.

\begin{fact}
If $(G, A)$ has relative property (T) and $\pi \: G \to G'$ is a homomorphism then $(\pi(G), \pi(A))$ has relative property (T). 
\end{fact}

\subsection{A Special Case}

We begin with the following lemma, which shows $(2) \implies (1)$ in the case when $A= \Z[S^{-1}]^n$.

\begin{lemma}
Suppose that $\G$ is a group and $ \f \: \G \to \GL{n}{\Z[S^{-1}]}$ a homomorphism such that $(\sdp{\G}{_\f \Z[S^{-1}]^n}, \Z[S^{-1}]^n)$ has relative property (T). Then $\~{\f(\G)}^Z(\Re)$ is non-amenable. 
\end{lemma}

\begin{proof}

Let $A = \Z[S^{-1}]^n$. Since $\ker(\f) \leq \sdp{\G}{A}$ centralizes $A$, it follows that $\ker(\f) \norm \sdp{\G}{A}$ and hence by Fact 2 $(\sdp{\f(\G)}{A}, A)$ has relative property (T).  

Recall that $A\leq V := \Re^n \times \Prod{p \in S} \Q_p^n$ is a co-compact lattice. So $(\sdp{\f(\G)}{V}, V)$ also has relative property (T) by Lemma 5.3.

Now since $\Prod{p \in S} \Q_p^n \norm \sdp{\f(\G)}{V}$ by Fact 2 we get that $(\sdp{\f(\G)}{\Re^n}, \Re^n)$ has relative property (T).  

This implies that $(\sdp{\~{\f(\G)}^Z(\Re)}{\Re^n}, \Re^n)$ has relative property (T). Indeed, any strongly continuous unitary representation of $\sdp{\~{\f(\G)}^Z(\Re)}{\Re^n}$ is a strongly continuous representation of $\sdp{\f(\G)}{\Re^n}$ (since $\f(\G)$ has the discrete topology).

But this means that $\~{\f(\G)}^Z(\Re)$ is non-amenable as is demonstrated by the next lemma.
\end{proof}

\begin{lemma}
Suppose that $G$ and $A$ are locally compact amenable groups such that $G$ acts on $A$ by automorphisms. Then $(\sdp{G}{A},A)$ has relative property (T) if and only if $A$ is compact. 
\end{lemma}

\begin{proof}

If $A$ is compact, then it has property (T) and hence $(\sdp{G}{A},A)$ has relative property (T).

Conversely, suppose that $(\sdp{G}{A}, A)$ has relative property (T). Recall that if $\mu_G$ and $\mu_A$ are (right invariant) Haar measures on $G$ and $A$ respectively, then $\mu_G \mu_A$ is a (right invariant) Haar measure on $\sdp{G}{A}$ (use Fubini's Theorem). Also recall that since $G$ and $A$ are amenable, the right regular representation $\rho \: \sdp{G}{A} \to \Un(L^2(\sdp{G}{A}))$ has almost invariant vectors. Then there is an $f \in L^2(\sdp{G}{A}) \bs \{0\}$ which is $A$ invariant, namely it is constant on the left cosets of $A$ and therefore it is a function of $G$ only. Then by Fubini's Theorem, 

\begin{equation}
\8  >  \int_{G \times A}|f(g)|^2 d\mu_G(g)d\mu_A(a) 
=\mu_A(A)\int_{G}|f(g)|^2 d\mu_G(g) > 0
\nonumber
\end{equation}

And therefore, $\mu_A(A) <\8$. This of course means that $A$ is compact.
\end{proof}


\subsection{The proof of Theorem 1' in the direction of $(2) \implies (1)$.}

\begin{proof}
Let $A$ be an Abelian group such that 
\begin{enumerate}
  \item The $\Q$-rank of $A$ is finite and non-zero.
  \item There is an action of $\G$ on $A$ by automorphisms such that $(\sdp{\G}{A},A)$ has relative property (T).
  \item The $\Q$-rank of $A$ is minimal among all Abelian groups satisfying (1) and (2). 
\end{enumerate}

Let $\mathrm{tor}(A) = \set{a \in A}{na = 0 \text{ for some } n \in \Z}$ be the torsion $\Z$-submodule of $A$. Observe that it is $\G$-invariant and hence $\mathrm{tor}(A) \norm \sdp{\G}{A}$. By Fact 2, we may assume that $A$ is torsion free. Since $\mathrm{tor}(A)$ is the kernel of the homomorphism $A \to \Q \otimes_\Z A$, we identify $A$ with it's image in $\Q \otimes_\Z A$. 

If $n$ is the $\Q$-rank of $A$ then there exists $v_1, \dots, v_n \in A$ such that $\Oplus{i = 1}{n}\Q\cdot v_i = \Q \otimes_\Z A$. (The notation is meant to emphasize the basis.) 

Now let $\f \: \G \to \GL{n}{\Q}$ be the corresponding homomorphism. (Observe that since $\G$ acts by automorphisms on $A \leq \Q \otimes_\Z A$ as an Abelian group, it acts by automorphisms of $A$ as a $\Z$-module. This means that we may extend the action $\Q$-linearly to obtain an automorphism of all $\Q \otimes_\Z A$.  And the group of automorphisms of $\Q \otimes_\Z A$, with respect to the above basis, is of course $\GL{n}{\Q}$.) 

Now, since $\G$ is finitely generated it follows by Lemma 8.1 that $\sdp{\G}{A}$ is finitely generated, and therefore $\sdp{\f(\G)}{A}$ is also finitely generated. So there is a finite set of primes $S_0$ such that $A \leq \Oplus{i = 1}{n}\Z[S_0^{-1}]\cdot v_i$.

For each $i = 1, \dots, n$, let $S_i = \set{p \in S_0}{ A \cap (\Z[S_0^{-1}] \cdot v_i )\subset \Q_p\cdot v_i \text{ is not precompact}}$. 

\begin{claim}
There is a $T \in \GL{n}{\Q}$ such that $T(A) \leq \Oplus{i = 1}{n}\Z[S_i^{-1}]\cdot v_i$ and $p \in S_i$ if and only if $T(A)\cap (\Z[S_i^{-1}]\cdot v_i) \subset \Q_p\cdot v_i$ is not precompact.
\end{claim}

\begin{proof}
For each $i = 1, \dots, n$ and $p \in S_0 \bs S_i$ there is a $k \in \N$ such that 
\begin{equation}
A\cap (\Z[S_0^{-1}]\cdot v_i) \subset \frac{1}{p^{k}}\Z_p\cdot v_i.
\nonumber
\end{equation} 
Let $k_i(p) \geq 0$ be the minimal one. Then, define the diagonal matrix:

\begin{equation}
T =  \(
\begin{array}{ccc}
    \Prod{p \in S_0 \bs S_1}p^{k_1(p)}  & &  0 \\
       & \ddots & \\
       0 & & \Prod{p \in S_0 \bs S_n}p^{k_n(p)}  \\
\end{array}
\)
\nonumber
\end{equation}

where of course we define $ \Prod{p \in S_0 \bs S_i}p^{k_i(p)} = 1$ in case $S_i = S_0$. 

Then, $T(A) \leq \Oplus{i = 1}{n}\Z[S_i^{-1}] \cdot v_i$ and $p \in S_i$ if and only if $T(A)\cap (\Z[S_i^{-1}]\cdot v_i) \subset \Q_p\cdot v_i$ is not precompact.
\end{proof}

Therefore, up to replacing $A$ by an isomorphic copy (and conjugating the $\G$-action), we may assume that $A \leq \Oplus{i = 1}{n}\Z[S_i^{-1}]\cdot v_i$ and that $p \in S_i$ if and only if $A\cap (\Z[S_i^{-1}]\cdot v_i )\subset \Q_p\cdot v_i$ is not precompact.

\begin{claim}
$A = \Oplus{i = 1}{n}\Z[S_i^{-1}]\cdot v_i$. 
\end{claim}

\begin{proof}
Let $i \in \{1, \dots, n\}$. Consider the set $C_i = \set{c \in \Z[S_i^{-1}]}{ cv_i \in A}$ which is a group under addition. Observe that $1 \in C_i$. 

We aim to show that $C_i =\Z[S_i^{-1}]$ and begin by showing that $\Z[\frac{1}{p}] \subset C_i$ for each $p \in S_i$. 

By definition, if $p \in S_i$ then for each $k \in \N$ there is a $c \in C_i$ such that $c = \frac{a}{bp^k}$ where $p$ does not divide $a$ and $b$. This means that $\frac{a}{p^k} = bc \in C_i$. Now, since $p$ does not divide $a$ it follows that there exists $x, y \in \Z$ such that $xp^k + ya= 1$. Namely, $x + y \frac{a}{p^k} = \frac{1}{p^k} \in C_i$. 

By induction, suppose that if $P \subset S$ is any subset of size $l-1$ that $\Z[P^{-1}] \subset C_i$. Then, for $p_1, \dots, p_l \in S_i$ and $k_1, \dots, k_l \in \N$ we have that 
\begin{equation}
\frac{1}{p_1^{k_1} \cdots p_{l-1}^{k_{l-1}}}, \frac{1}{p_2^{k_2} \cdots p_{l}^{k_{l}} } \in C_i
\nonumber
\end{equation}
Since $p_1$ and $p_l$ are relatively prime, there exists $x, y \in \Z$ such that $xp_l^{k_l} + yp_1^{k_1} = 1$. Then, 
\begin{equation}
\frac{x}{p_1^{k_1} \cdots p_{l-1}^{k_{l-1}}} + \frac{y}{p_2^{k_2} \cdots p_{l}^{k_{l}}} = \frac{xp_l^{k_l} + yp_1^{k_1} }{p_1^{k_1} \cdots p_{l}^{k_{l}}} = \frac{1}{p_1^{k_1} \cdots p_{l}^{k_{l}}} \in C_i
\nonumber
\end{equation}
\end{proof}

Observe that this means that for an arbitrary $v = \Sum{i = 1}{n} \a_i v_i \in \Q\otimes_\Z A$ we have that $v \in A$ if and only if $\a_i \in \Z[S_i^{-1}]$ for each $i = 1, \dots, n$.

Now, up to renumbering the basis, assume that $|S_1| \geq |S_i|$ for each $i = 1, \dots, n$ and \\$S_1 = \cdots = S_m$ and $S_1 \neq S_i$ for any $ i = m+1, \dots, n$. Let $S = S_1$.

\begin{claim}
The subgroup $\Oplus{i = 1}{m}\Z[S^{-1}]\cdot v_i$ is $\G$-invariant.
\end{claim}

\begin{proof}
Let $\g = \( \g_{i,j}\)$ be the matrix representation of $\g$ with respect to the above basis. Observe that $\Oplus{i = 1}{m}\Z[S^{-1}]\cdot v_i$ is $\G$-invariant if and only if for every $(\g_{i,j}) \in \G$ and each $i_0 \in \{1, \dots, m\}$

\begin{equation}
\nonumber
\g_{j_0,i_0} \in
\begin{cases}
    \Z[S^{-1}]  & \text{if } j_0 \in\{ 1, \dots, m\}, \\
    \{0\}  & \text{if } j_0 \in \{m+1,\dots, n\}.
\end{cases}
\end{equation}

Since $\G$ preserves $A$ the above condition is already satisfied for $j_0 \in \{1, \dots, m\}$. We now show that if $i_0 \in \{ 1, \dots, m\}$ and $j_0 \in \{ m+1, \dots, n\}$ then $\g_{j_0,i_0} = 0$.

By maximality of $|S|$ and the fact that $S \neq S_{j_0}$ there is a $p \in S \bs S_{j_0}$. Now, $\frac{1}{p^l}v_{i_0} \in A$ for each $l \in \N$ so that $\g(\frac{1}{p^l}v_{i_0}) \in A$ as well. 

This means that $\frac{1}{p^l}\g_{j_0,i_0} \in \Z[S_{j_0}^{-1}]$ and so $\frac{m}{p^l}\g_{j_0,i_0} \in \Z[S_{j_0}^{-1}]$  for every $m \in \Z$. 

Choose $m \in \Z \bs \{0\}$ and $l \in \N$ sufficiently large such that 
\begin{equation}
\nonumber
\frac{m}{p^l} \g_{j_0,i_0} \in \Z[S_{j_0}^{-1}] \cap \Z[p^{-1}] = \{0\}.
\end{equation}

\end{proof}

We are almost done. Indeed the result follows by Lemma 8.2 and the following:

\begin{claim}
Let $A'= \Oplus{i = 1}{m}\Z[S^{-1}]\cdot v_i$. Then $A' =A$.
\end{claim}

\begin{proof}
If we can show that the $\Q$-rank of $A/A' \cong \Oplus{i = m+1}{n}\Z[S_i^{-1}]\cdot v_i$ is 0 then the result follows.

Since $A'$ is $\G$-invariant it follows that $A' \norm \sdp{\G}{A}$. By Fact 2, $(\sdp{\G}{(A/A')}, A/A')$ has relative property (T). However, $A$ was chosen to be of minimal (non-zero) $\Q$-rank among all such Abelian groups and so the $\Q$-rank of $A/A'$ is 0.
\renewcommand{\qedsymbol}{}
\end{proof}

\end{proof}


\section{Some Examples}

We would like to take the opportunity to address two questions that may naturally arise as one reads this exposition.

\begin{questions}
Does every nonamenable linear group satisfy condition (1) of Theorem 1? Namely, if $\G$ is a non-amenable linear group does there always exist $\f \: \G \to \SL{n}{\Re}$ with $\~{\f(\G)}^Z(\Re)$ non-amenable? 
\end{questions}

The answer to this question is of course no. There are purely $p$-adic higher rank lattices and by Margulis' Superrigidity Theorem such lattices only admit precompact homomorphisms into $\SL{n}{\Re}$ \cite[Example IX (1.7.vii) p. 297, Theorem VII (5.6)]{Margulis}.

The second question arises out of the following application:

\begin{theorem}[\cite{PopaGab}, \cite{Tornquist}, \cite{Shalom}{p23}]
Let $\a \: \G \to \Aut(A)$ be a homomorphism, with $A$ discrete Abelian such that $(\G \ltimes _\a A, A)$ has relative property (T). Then there are uncountably many orbit inequivalent free actions of the free product $\a(\G) \star \Z$ on the standard probability space. 
\end{theorem}

We point out that although both the papers of Gaboriau-Popa and T\"ornquist  prove the above theorem for the case of $A = \Z^2$ and $\G = F_n$, Y. Shalom points out that the proof extends to show the above theorem.

Theorem 9.1, taken with Theorem 1, shows that it is good to know if such semidirect products may be constructed with the action of $\G$ on the Abelian group $A$ being faithful.

\begin{questions}
Does there exist a linear group $\G$ satisfying property (F$_\8$) such that every homomorphism $\f \: \G \to \SL{n}{\Q}$ is not injective?
\end{questions}

The answer to this question is yes. The homomorphism $\f'$ found in the proof of Theorem 1 will have a kernel in general. This kernel arises out of the need to specialize transcendental extensions of $\Q$ in order to get an action on an Abelian group of finite $\Q$-rank. We therefore look to these transcendental extensions to find our example.

\begin{prop}
Every homomorphism $\f \: \SL{3}{\Z[x]} \to \GL{n}{\Q}$ is not injective.
\end{prop}

We remark that this proposition only shows that $\SL{3}{\Z[x]}$ never has a faithful action on an Abelian group of finite $\Q$-rank. On the otherhand, it is possible to get relative proeprty (T) from this group. Indeed, Y. Shalom showed \cite[Theorem 3.1]{Shalom1999} that $(\SL{3}{\Z[x]} \ltimes \Z[x]^3, \Z[x]^3)$ has relative property (T).

To prove this proposition, we will need the following:

\begin{definition}
Let $\G$ be a group generated by the finite set $S$. An element $\g \in \G$ is said to be a $U$-element if 
\begin{equation}
d_S(\g^m, 1) = O(\log m)
\nonumber
\end{equation}
 where $d_S$ is the metric on the $S$-Cayley graph of $\G$ and 1 is of course the identity.
\end{definition}
 
This property is wonderful because it identifies ``unipotent'' elements while appealing only to the internal group structure. This is exemplified by the following:

\begin{prop}[\cite{LMR}{Proposition 2.4}] If $\g\in \G$ is a $U$-element then for every representation $\f \: \G \to \GL{n}{\Re}$ we have that $\f(\g)$ is virtually unipotent. 
\end{prop}

We now turn to the proof of Proposition 9.2:

\begin{proof}

Let $E_{i,j}(y)$ be the elemtary unipotent matrix in $\SL{3}{\Z[x]}$ with $y \in \Z[x]$ in the $(i,j)$-th position, and $i \neq j$. It is by now a well known result of Bass, Milnor and Serre (\cite[Corollary 4.3]{BMS}) that $\SL{3}{\Z}$ is generated by $S_1 := \{E_{i,j}(1)\}$. A similar result of Suslin (\cite{Suslin}) states that $\set{E_{i,j}(y)}{y \in \Z[x]}$ generates $\SL{3}{\Z[x]}$. By observing that, for a fixed $y \in \Z[x]$, all the $E_{i,j}(y)$ are conjugate (in $\SL{3}{\Z}$), and the following commutator relation, we see that the finite set $S_x := \{E_{i,j}(x)\}\cup S_1$ actually generates $\SL{3}{\Z[x]}$:

\begin{equation}
[E_{1,2}(y_1), E_{2,3}(y_2)] = E_{1,3}(y_1 y_2). 
\nonumber
\end{equation}

\begin{claim}
$E_{1,3}(y)$ is a $U$-element for each $y \in \Z[x]$.
\end{claim}

\begin{proof}
By Corollary 3.8 of \cite{LMR} $E_{i,j}(1)$ is a $U$-element. Furthermore, observe that  $d_{S_x}(E_{1,2}(m), 1) \leq d_{S_1}(E_{1,2}(m), 1)$ since $S_1 \subset S_x$. For $m$ sufficiently large, the above commutator relation, with $y_1 = m$ and $y_2 =y$, gives us that 
\begin{equation}
d_{S_x}(E_{1,3}(my), 1) \leq 2d_{S_x}(E_{1,2}(y), 1)+ 2d_{S_x}(E_{1,2}(m), 1) \leq 2(1+C )\log m
\nonumber
\end{equation}

where $d_{S_x}(E_{1,2}(m), 1) \leq C \log m$. Hence $E_{1,3}(y)$ is a $U$-element. 
\end{proof} 

Now to conserve notation, for each $y \in \Z[x]$ let us define $\g_y = E_{1,3}(a_y y)$ where $a_y \in \N$ is the minimum of all $a \in \N$ such that $\f(E_{1,3}(ay))$ is unipotent.

Also, let $G_u := \~{\<\f( \g_y) | y \in \Z[x] \>}^Z$ be the Zariski-closure. Then $G_u$ is $\Q$-rationally isomorphic to $\Re^d$ for some $d$. Indeed, $G_u$ is a $\Q$-group generated by commuting unipotent elements and is therefore both unipotent and Abelian. This means that there is a $\Q$-basis of $\Re^n$ for which $G_u$ is a subgroup of the upper triangular unipotent matrices, which is in turn isomorphic to $\Re^{n-1} \ltimes \Re^{n-2} \ltimes \cdots \ltimes \Re$. 

Now, fix a $\Q$-rational isomorphism $\rho \: G_u \to \Re^d$. Then, since $\set{\rho\phi(\g_y)}{ y \in \Z[x]}$ is Zariski-dense in $\Re^d$ there exists $y_1, \dots, y_d \in \Z[x]$ so that $\{\rho\phi(\g_{y_1}), \dots, \rho\phi(\g_{y_d})\}$ is a $\Q$-basis for $\Re^d$. 

Let $y \in \Z[x]$ such that $\<y\> \cap \set{\Sum{j = 1}{d}a_j y_j}{a_j \in \Z} = \{0\}$. Since $\rho\phi(\g_y)$ is in the $\Q$-span of our basis, there exists $q_j \in \Q$ such that
\begin{equation}
\rho \phi (\g_y) = \Sum{j = 1}{d} q_j \rho \phi (\g_{y_j}).
\nonumber
\end{equation}

Clearing the denominators we have that there are $m, m_1, \dots, m_d \in \Z$ such that

\begin{equation}
\g_y^m \Prod{j = 1, \dots, d}\g_{y_j}^{m_j} = E_{1,3}\(ma_y y + \Sum{j = 1}{d} m_j a_{y_j} y_j\) \in \ker(\rho \circ \phi).
\nonumber
\end{equation}
By our choice of $y$ and the fact that $\ker(\rho) =1$, we have that $\ker(\f) \neq1$.
\end{proof}

\bibliography{Rel(T)bib}
\bibliographystyle{alpha}

\end{document}